\documentclass[aap,MSNbibl,seceqn,dvips]{arximspdf}
\usepackage{graphicx}
\doi{10.1214/12-AAP876} 
\volume{23}
\issue{4}
\pubyear{2013}
\firstpage{1437}
\lastpage{1468}

\makeatletter
\newcommand{\rrvert}{\vert}
\newcommand{\llvert}{\vert}
\newcommand{\eqref}[1]{(\ref{#1})}
\newtheorem{theorem}{Theorem}
\newtheorem{corollary}[theorem]{Corollary}
\newtheorem{lemma}{Lemma}
\newtheorem{proposition}{Proposition}
\newcommand{\ve}{\varepsilon}
\newcommand{\R}{\mathbb{R}}
\newcommand{\EE}{\mathbb{E}}
\newcommand{\Prob}{\mathbb{P}}
\newcommand{\bphi}{\bolds{\phi}}

\makeatother

\begin{document}
\begin{frontmatter}

\title{Dynamics of cancer recurrence\thanksref{TT1}}
\thankstext{TT1}{Supported by the US National Science Foundation under Grant NSF DMS-12-24362.}
\runtitle{Dynamics of cancer recurrence}

\begin{aug}
\author[A]{\fnms{Jasmine} \snm{Foo}\corref{}\ead[label=e1]{jyfoo@math.umn.edu}}
\and
\author[B]{\fnms{Kevin} \snm{Leder}\ead[label=e2]{kevin.leder@isye.umn.edu}}
\runauthor{J. Foo and K. Leder}
\affiliation{University of Minnesota}

\address[A]{School of Mathematics\\
University of Minnesota\\
206 Church St SE\\
Minneapolis, Minnesota 55455\\
USA\\
\printead{e1}}

\address[B]{Department of Industrial\\
\quad and Systems Engineering\\
University of Minnesota\\
111 Church St SE\\
Minneapolis, Minnesota 55455\\
USA\\
\printead{e2}}

\end{aug}

\received{\smonth{1} \syear{2012}}
\revised{\smonth{5} \syear{2012}}

%
\begin{abstract}
Mutation-induced drug resistance in cancer often causes the failure of
therapies and cancer recurrence, despite an initial tumor reduction.
The timing of such cancer recurrence is governed by a balance between
several factors such as initial tumor size, mutation rates and growth
kinetics of drug-sensitive and resistance cells. To study this
phenomenon we characterize the dynamics of escape from extinction of a
subcritical branching process, where the establishment of a clone of
escape mutants can lead to total population growth after the initial
decline. We derive uniform in-time approximations for the paths of the
escape process and its components, in the limit as the initial
population size tends to infinity and the mutation rate tends to zero.
In addition, two stochastic times important in cancer recurrence will
be characterized: (i) the time at which the total population size first
begins to rebound (i.e., become supercritical) during treatment, and
(ii) the first time at which the resistant cell population begins to
dominate the tumor.
\end{abstract}

%


\end{frontmatter}

\section{Introduction}\label{sec1}

We consider a situation arising from population genetics, where a
population with net negative growth rate can escape certain extinction
via creation of a new mutant type. This scenario arises in a variety of
biological and medical applications. In particular, we consider the
following scenario in which a population of drug-sensitive cancer cells
is placed under therapy, leading to a sustained overall decline in
tumor size. However, drug-resistance mutations may arise in the
population, conferring a net positive growth rate to mutated cells and
their progeny under therapy. If a mutant arises prior to extinction of
the original population and forms a viable, growing subpopulation, then
the population has ``escaped'' extinction. These types of escape events
due to acquired resistance cause the failure of many drugs including
antibiotics, cancer therapies and anti-viral therapies. In the cancer
setting, the discovery of new molecularly targeted therapies has lead
to dramatic successes in tumor reduction in the past decade; however,
the majority of these therapies fail due to the development of drug
resistance and subsequent increase in tumor burden and progression of
disease. Examples of targeted therapies for which acquired resistance
exists include erlotinib/gefitinib in EGFR-mutant nonsmall cell lung
cancer, imatinib, dasatinib or nilotinib in BCR--ABL driven chronic
myeloid leukemia and vemurafenib in BRAF-mutant melanoma.

There has been a significant amount of previous work in the cancer
modeling literature on understanding the evolutionary dynamics of drug
resistance in cancer. For example, using stochastic processes with a
differentiation hierarchy to represent the sensitive and resistant
cells of a tumor, Coldman and Goldie studied the emergence of
resistance to one or two drugs~\cite{CoGo86,GoCo79,GoCo83}. In a
different twist, Harnevo and Agur studied drug resistance emerging due
to oncogene amplification using a stochastic branching process model
\cite{HaAg92,HaAg91}. Others have used multi-type branching process
models to study the probability of resistance emerging due to point
mutations in a variety of situations, for example,~\cite{Mi06,IwMiNo03}.
Komarova and Wodarz also utilized a multi-type branching
model to investigate the general situation in which~$k$ mutations are
required to confer resistance against $k$ drugs~\cite{Ko06,KoWo05}.
Most recently, in~\cite{FoMi2010} the authors considered an
inhomogeneous process wherein the birth and death rates of both
sensitive and resistant cells are dependent upon a temporally varying
drug concentration profile, to accommodate the effects of
pharmacokinetic dynamics as the drug is metabolized over time. The
analysis in most of these works has been focused on calculations of the
eventual probability of developing resistance and the resistant
population size, rather than the variable timing of tumor recurrence.

In addition to work specifically related to mathematical modeling of
cancer recurrence, we also discuss some mathematical contributions to
the study of extinction paths in subcritical branching processes and
the dynamics of escape in this context. In particular, in \cite
{JaKlSa07Adv} Jagers and co-authors considered large population
approximations of ``the path to extinction'' in Markovian sub-critical
branching processes. In this work they established convergence of
finite dimensional distributions of these paths viewed on the time
scale of extinction. The follow-up work~\cite{JaKlSa07PN} generalized
these results to a broader class of inter-arrival times (i.e.,
distributions more general than exponential). Sehl et al. investigated
the limiting moments of extinction times of subcritical branching
processes, and used this as a tool for investigating the effects of
various cancer therapies on healthy tissue~\cite{SeZhSiLa11}. Last,
Sagitov and Serra characterized the asymptotic structure for BGW
process with escape, as mutation rate $\mu\rightarrow0$, conditioned
on successful escape, which is an important asymptotic regime in many
problems such as the evolution of new species
\cite{SaSe09}.

A typical solid tumor has a density between $10^7$ and $10^9$ cancer
cells per cubic centimeter~\cite{MiHaKo99}. Therefore, in this work we
are interested in deriving path approximations of the escape process
that are uniform in time, in the regime of a very large initial
population. In the large population limit, it is tempting to assume
that the stochastic model can be approximated by a purely deterministic
model. However, a simple comparison of the mean behavior of the
stochastic model with a deterministic model illustrates that it is
important to consider the stochasticity of the extinction time. Here,
we develop limiting stochastic approximations for the population
process that greatly simplify the population process model while
maintaining the stochastic extinction time behavior. Interesting
earlier work by Jagers, Sagitov et al. established convergence of the
finite dimensional distributions for the declining sensitive cell
populations on the time scale of extinction, leaving open the question
of tightness~\cite{JaKlSa07Adv}. In the present work we first construct
nearly-deterministic uniform in time limit approximations to both the
declining sensitive population paths as well as the supercritical
resistant cell escape paths. Then, tightness of the joint sensitive and
resistant process can be established as a simple consequence of these
approximations, yielding the weak convergence result in simpler fashion
than via direct analysis of the joint process.\looseness=-1

\begin{figure}

\includegraphics{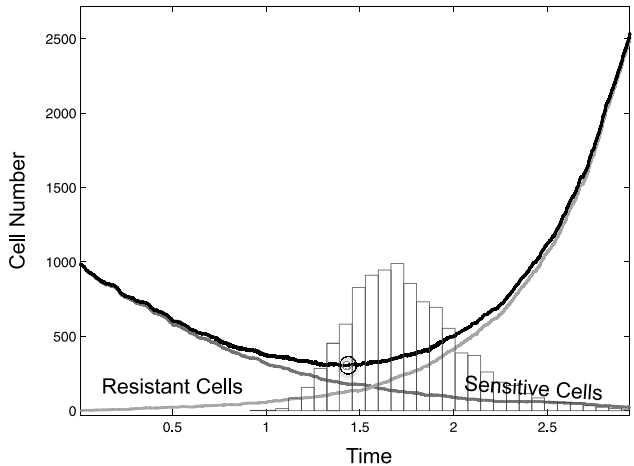}

\caption{Sample simulation of escape dynamics (population size
versus time). Dark black line: total population size of the tumor,
labeled grey lines: resistant and sensitive cell population size. The
circle marks the minimum of the total tumor size process (i.e., the
turnaround time), and the point at which red and blue lines cross is
the crossover time. A histogram plotting the distribution of the
minimum turnaround time for each sample path is plotted in green in
the background. Parameters: starting population 1000 sensitive, 0
resistant. Net growth rate of sensitive and resistant birth--death
processes are $-1.0$ and 2.0, respectively, and the mutation rate $\mu=
0.01$.}\label{figsamplesim}
\end{figure}

We then use these approximations to characterize the distribution of
``turnaround times,'' at which the total population size switches from
subcritical to supercritical. In the clinical context, this represents
the time at which progression of disease is observed through serial
tumor scans or bloodwork (in leukemias); thus the ability to
characterize and predict this time is of significant prognostic
interest. In addition, we characterize the ``crossover time'' at which
the resistant mutants first overtake the original type in the
population. Estimates of crossover times, and more generally the times
at which certain composition thresholds are reached, are extremely
useful in clinical decision-making. For example, when simultaneous
combination therapies are considered, understanding these random times
allows for informed decisions on the optimal time to switch to another
therapy and thus ``target'' a different subpopulation of cells within
the tumor. Figure~\ref{figsamplesim} illustrates these times in a
sample path simulation of the process, in addition to a sample
distribution of turnaround times. Our results are derived in the
framework where the time scale of the processes is sped up by the
extinction time of the original population, a natural time scale since
this time represents the maximum length of effectiveness of the drug.
We restrict our attention to binary branching processes which are
appropriate for modeling cancer cell populations undergoing binary
division; however, these results can be extended to study more general
offspring distributions, and thus may be useful for studying escape
dynamics in viral populations, for instance.

The rest of the paper is organized as follows. In Section~\ref{sec2} we
introduce the model and discuss earlier results in the field. In
Section~\ref{sec3} we present some results on the mean of the resistant cell
population at multiples of the extinction time. In Section~\ref{sec4} we present
a path approximation result where we show that the limit process
uniformly approximates both the sensitive and resistant cell process on
the time scale of the extinction time of the sensitive cells. We
determine limiting distributions of the crossover time when the
resistant cell population first becomes dominant, and the random time
of disease progression or the ``turnaround'' time. In Section~\ref{sec5} we
briefly illustrate an application of these results to studying the time
of disease recurrence due to drug resistance in nonsmall cell lung
cancer (NSCLC). In Section~\ref{secproofs} we present the proofs of our main results.

Throughout the paper we use the following standard Landau asymptotic
notation for nonnegative functions $f (
\cdot ) $ and $g ( \cdot )$: $f ( x ) =O (
g (
x )  ) $ means that $f ( x ) \leq cg ( x ) $
for some $c\in ( 0,\infty ) $, $f ( x ) =\Omega (
g ( x )  ) $ if and only if $f ( x ) \geq
cg ( x ) $, $f ( x ) =o ( g ( x )
) $ holds if and only if $f ( x ) /g ( x )
\rightarrow0$ as $x\rightarrow\infty$ and last, $f(x)\sim g(x)$ holds
if and only if $f(x)/g(x)\to1$ as $x\to\infty$.

\section{Model and previous work}\label{sec2}
In this section we introduce the mathematical model and notation, and
review previous results on related problems.
We start with an initial population of drug sensitive cells with size
$x$. This population $Z_0(t)$ is modeled as a subcritical Markovian
binary branching process which declines during treatment with net
growth rate $\lambda_0 < 0$, birth rate $r_0$ and death rate $d_0$; we
will also use the notation $|\lambda_0|=r$. Resistance mutations arise
at rate $\mu_xZ_0(t)$, and each of these mutations gives rise to a
supercritical Markovian binary branching process initialized by one
mutant cell with net growth rate $\lambda_1>0$. We set $\mu_x=\mu
x^{-\alpha}$ for $\mu>0$ and $\alpha\in(0,1)$. The total population of
mutants, which we will call ``resistant cells,'' is denoted $Z_1(t)$.
These processes are defined on a probability space $(\Omega,\mathcal
{F},\Prob)$. In addition, we define the filtration $\mathcal{F}_t^i$
generated by $Z_j(s)$, for $s\leq t$ and $j\leq i$. Note that in this
work, unless otherwise stated, the expectation and probability
operators are conditioned on the initial conditions $Z_0(0) = x$ and
$Z_1(0)= 0$.

Since the net growth rate of the original population is negative, it
will go extinct eventually with probability 1. We will denote this time
of extinction by $T_x$, where $x$ denotes the starting population. The
following limit theorem from~\cite{Pa89} will prove useful throughout
the rest of the paper:
%
\begin{equation}
\label{eqexttimeconv} T_x-\frac{1}{r}\log x\Rightarrow
\frac{1}{r} (\eta+\log c )\qquad\mbox{as }x\to\infty,
\end{equation}
where $\eta$ is a standard Gumbel random variable and $c$ is the Yaglom
constant for~$Z_0$. For a binary branching process, the Yaglom constant
has the form $(d_0-r_0)/r_0$.

Previously, Jagers and colleagues~\cite{JaKlSa07Adv} studied the paths
to extinction in a subcritical Markovian branching process, which we
will also call~$Z_0$ starting at size~$x$. They considered the process~$Z_0$
on the time scale of the extinction time and established
convergence in finite dimensional distributions as \mbox{$x \rightarrow\infty$}.
%
\begin{theorem}[(Jagers et al.~\cite{JaKlSa07Adv})]
For $u \in[0,1)$,
\[
x^{u-1}Z_0(uT_x) \mathop{\rightarrow}^{FD}
c^{-u} e^{-u\eta}.
\]
\end{theorem}
Similar results on convergence in finite dimensional distribution of
subcritical branching processes with more general inter-arrival times
were also shown in~\cite{JaKlSa07PN}. In addition, Kimmel and Wu
generalized these results to consider the case of critical branching
processes~\cite{WuKi10}.

\section{Mean of $Z_1(uT_x)$}\label{sec3}
In this section we examine the growth rate of the mean of $Z_1$. In
addition, we examine a common modeling assumption and note the
importance of considering the tails of the extinction time $T_x$ in
studies of escape dynamics.
We will first consider the expected resistant population at $vT_x$ for
some $v>0$ (and temporarily assume $\alpha=0$),
%
\begin{equation}
\label{eqexpmuts} \EE \bigl[Z_1(vT_x) \bigr]= \EE
\biggl[\mu T_x\int_0^{v\wedge
1}Z_0(uT_x)
\exp \bigl(\lambda_1T_x(v-u) \bigr)\,du \biggr].
\end{equation}
If we assume that sensitive cells follow a deterministic decay
$Z_0(t)=xe^{\lambda_0 t}$ and approximate their extinction time as
$T_x\approx-\frac{1}{\lambda_0}\log x$, then we can heuristically
estimate the expected value as
\begin{eqnarray*}
\EE\bigl[Z_1(vT_x)\bigr] &=& \frac{\mu}{r}\log x
\int_0^{v\wedge1}x^{1-u}x^{({\lambda_1}/{r})(v-u)}\,du
\\
&=& \frac{\mu}{r}x^{1-{\lambda_1}v/{\lambda_0}}\log x\int_0^{v\wedge
1}x^{-u(1+{\lambda_1}/{r})}\,du
\\
&=& \frac{\mu}{\lambda_1-\lambda_0}x^{1+{\lambda_1}v/{r}} \biggl(1-\exp \biggl[-(v\wedge1) \biggl(1+
\frac{\lambda_1}{r}\biggr)\log x \biggr] \biggr).
\end{eqnarray*}
Thus we observe that this expected value is finite for all $v>0$.

However, suppose that there is just a single sensitive cell at time
$t=0$ whose birth rate is $r_0=0$, and death rate $r$. Then, of
course, the extinction time satisfies $T_1\sim\exp(r)$ and by
conditioning on this time we have
\[
\EE\bigl[Z_1(vT_x)\bigr]=\mu\int_0^{\infty}
\EE\bigl[e^{\lambda_1(vT_1-s)}1_{\{T_1\geq
v/s\}}\bigr]\,ds.
\]
Due to the exponential tails of $T_1$ we see that the above integral
diverges to $\infty$ for $\lambda_1v\geq r$ and we clearly see the
importance of the randomness in the extinction time.
%
%
The previous result easily applies to models with births and deaths in
the sensitive cell population. In particular, we have the following proposition.
%
\begin{proposition}\label{propinfmean}
Let $v>0$ and $\frac{\lambda_1 v}{r} > 1$, then for all $x$, $\EE
[Z_1(vT_x)]$ is infinite.
\end{proposition}


\begin{pf}
By conditioning on $Z_0(s), s>0$ and then applying a change of measure
we can write the integral of interest as
\begin{eqnarray*}
\EE \bigl[Z_1(vT_x) \bigr] &=& \EE \biggl[\mu
T_x\int_0^{v\wedge
1}Z_0(uT_x)
\exp \bigl[\lambda_1T_x(v-u) \bigr]\,du \biggr]
\\
&=& \mu\int_0^{v\wedge1}\int_0^{\infty}te^{\lambda_1 t(v-u)}
\EE \bigl[Z_0(uT_x)| T_x\in dt
\bigr]g_x(t)\,dt \,du,
\end{eqnarray*}
where $g_x(t)\,dt=\Prob(T_x\in dt)$.

Noting the fact that if $u<1$ then $E [Z_0(uT_x)| T_x\in dt
] \geq1$, we can bound this from below by
\begin{eqnarray*}
\EE \bigl[Z_1(vT_x) \bigr] &\geq&\mu\int
_0^{v\wedge1}\int_0^{\infty
}te^{\lambda_1 t(v-u)}g_x(t)\,dt
\,du
\\
&=& \mu\int_0^{v\wedge1} \biggl( \EE[T_x] +
\lambda_1 (v-u) \int_0^{\infty}e^{\lambda_1 s(v-u)}
\int_s^\infty t g_x(t)\,dt \,ds \biggr)
\,du
\\
&\geq&\mu\int_0^{v\wedge1} \biggl( \EE[T_x]
+ \lambda_1 (v-u) \int_0^{\infty}e^{\lambda_1 s(v-u)}s
\int_s^\infty g_x(t)\,dt \,ds \biggr) \,du
\\
&\geq& \mu\int_0^{v\wedge1}
\biggl( \EE[T_x] + c\lambda_1 (v-u) \int
_{t_0}^{\infty}e^{\lambda_1 s(v-u)}se^{-rs} \,ds
\biggr) \,du.
\end{eqnarray*}
The final inequality is based on the fact that for $x\geq1$, $\Prob
(T_x>s) \geq\Prob(T_1>s)$ and the asymptotic result that as $t\to\infty
$, $\Prob(T_1>t)\sim ce^{-rt}$.
Considering the final equation in the previous display, we see that if
$\lambda_1v>r$ then for $u$ sufficiently small, the inner integral
diverges to $\infty.$
\end{pf}

We can easily find the asymptotic growth rate of $\EE[Z_1(vT_x)]$ as
$x\to\infty$. Based on the previous subsection, we know that this is
only meaningful if we consider $v\leq-\lambda_0/\lambda_1$; for
simplicity we will just assume that $v\leq1$ and $r=|\lambda_0|\geq
\lambda_1$. Earlier heuristic calculations (where we set $\alpha=0$)
indicate that the mean of $Z_1(vT_x)$ grows like $x^{1+v\lambda_1/r}$
as $x\nearrow\infty$.
In particular we have the following theorem.
%
\begin{theorem}
\label{thmmeanscaling}
Assume that $r\geq\lambda_1$, then for $v\in(0,1]$ and $\alpha\in
(0,1)$ we have that
\[
\EE \bigl[Z_1(vT_x) \bigr]\sim x^{1+v\lambda_1/r-\alpha}
\frac{c^{\lambda
_1v/r}\mu\Gamma(1-{\lambda_1v}/{r})}{ \lambda_1+r}.
\]
\end{theorem}
We defer the proof of this result to Section 7.

\section{Paths of escape}\label{sec4}
We now establish an approximation theorem for the paths of the joint process
$(Z_0(uT_x), Z_1(uT_x))$. In the large $x$ limit, scaled versions of
these paths can be approximated uniformly in time by a simple
stochastic process whose only source of randomness arises from the
stochasticity of the limit theorem for the extinction time.

Before beginning, we first establish some notation. We will work with
scaled versions of the sensitive and resistant populations sped up in
time. Let us define $s_x(t) = \frac{1}{r} \log x + t$. For $u\in[0,1]$
and $t\in\R$, define
%
\begin{eqnarray}
\label{eqscaleprocessdef} Z_0^x \bigl(us_x(t)
\bigr) &=& x^{u-1}Z_0 \biggl(u \biggl(\frac{1}{r}\log
x+t \biggr) \biggr),
\nonumber
\\[-8pt]
\\[-8pt]
\nonumber
Z_1^x \bigl(us_x(t) \bigr) &=&
x^{-\lambda_1u/r-1 + \alpha}Z_1 \biggl(u \biggl(\frac{1}{r}\log x+t
\biggr) \biggr).
\nonumber
\end{eqnarray}

Throughout the rest of the paper, the superscript $x$ will denote
scaling by the appropriate function of $x$. For ease of notation we
introduce the following notation:
\begin{eqnarray*}
\phi_0^x(u,t) &=& \EE Z_0^x
\bigl(us_x(t) \bigr) = e^{\lambda_0 ut},
\\
\phi_1^x(u,t) &=& \EE Z_1^x
\bigl(us_x(t) \bigr) = \frac{\mu e^{\lambda_1
ut}}{\lambda_1-\lambda_0} \bigl( 1 -
e^{(\lambda_0-\lambda_1)ut}x^{
{(\lambda_0-\lambda_1)u}/{r}} \bigr).
\end{eqnarray*}
In addition, we will sometimes need to work with the population
processes sped up in time but not scaled in space, which are defined
for ${Z}_i(us_x(t))$, for $i=0,1$
and their means are denoted by
\begin{eqnarray*}
\phi_i(u,t) = \EE Z_i\bigl(us_x(t)\bigr).
\end{eqnarray*}

In the following, we establish the approximation result by first
showing that for any $t\in\mathbb{R}$ we can approximate the scaled
joint process by its mean uniformly in~$u$. This is done by martingale
arguments and showing relevant second moments are uniformly bounded in
$x$. We then prove that this approximation is uniform for~$t$ in
compact sets, and that one can approximate $(Z_0(uT_x),Z_1(uT_x))$
uniformly in time by $(\phi_0(u,T_x-\frac{1}{r}\log x),\phi_1(u,T_x-\frac{1}{r}\log x))$, where the previous formula is
interpreted as the mean functions $\phi_i^x$ evaluated at the random
parameter $T_x-\frac{1}{r}\log x$. We begin with a result on the
moments of $Z_0$ and $Z_1$.
%
\begin{lemma}
\label{lemmamoments}
Let $\tilde{Z}_1$ be a binary branching process starting from size one
with birth rate $r_1$ and death rate $d_1$, then for $0<s<s_x(t)$,
\begin{longlist}[(iii)]
\item[(i)]
\begin{eqnarray*}
\EE\bigl[Z_1 \bigl(s_x(t) \bigr)^2\bigr]&=&
\frac{\mu^2}{x^{2\alpha}}\int_0^{s_x(t)}\int
_0^{s_x(t)}\EE\bigl[Z_0(s)Z_0(y)
\bigr]e^{\lambda
_1(s_x(t)-s)}e^{\lambda_1(s_x(t)-y)}\,ds\,dy
\\
&&{} + \frac{\mu}{x^{\alpha}}\int_0^{s_x(t)}\EE
Z_0(s)\EE\bigl[\tilde {Z}_1\bigl(s_x(t)-s
\bigr)^2\bigr]\,ds.
\end{eqnarray*}
\item[(ii)]
\[
\EE\bigl[Z_0(s)Z_1 \bigl(s_x(t)
\bigr)\bigr]=\frac{\mu}{x^{\alpha}}\int_0^{s_x(t)}\EE
\bigl[Z_0(y)Z_0(s)\bigr]e^{\lambda_1(s_x(t)-y)}\,dy.
\]
\item[(iii)]
\begin{eqnarray*}
\operatorname{Var}\bigl[Z_1\bigl(s_x(t)\bigr)\bigr]&=&
\frac{\mu^2}{x^{2\alpha}}\int_0^{s_x(t)}\int
_0^{s_x(t)}\operatorname{Cov}\bigl(Z_0(s),Z_0(y)
\bigr)e^{\lambda_1(2s_x(t)-(s+y))}\,ds\,dy
\\
&&{} + \frac{\mu}{x^{\alpha}}\int_0^{s_x(t)}\EE
Z_0(s)\EE\bigl[\tilde {Z}_1\bigl(s_x(t)-s
\bigr)^2\bigr]\,ds.
\end{eqnarray*}
\end{longlist}
\end{lemma}
The proof of this result can be found in Section~\ref{secproofs}.

Lemma~\ref{lemmamoments} allows us to establish the following result
via the Doob's maximal inequality.
%
\begin{lemma}
\label{lemmaapprox1}
For $a\in(0,1)$, $\ve>0$ and $t\in\mathbb{R}$,
\begin{longlist}[(ii)]
\item[(i)]
\[
\lim_{x\to\infty} \mathbb{P} \Bigl(\sup_{u\in[0,a]} \bigl\llvert
Z_0^x\bigl(us_x(t)\bigr)-
\phi_0^x(u,t) \bigr\rrvert >\ve \Bigr) = 0.\vadjust{\goodbreak}
\]
\item[(ii)]
\[
\lim_{x\to\infty} \mathbb{P} \Bigl(\sup_{u\in[0,1]} \bigl\llvert
Z_1^x\bigl(us_x(t)\bigr)-
\phi_1^x(u,t) \bigr\rrvert >\ve \Bigr) = 0.
\]
\end{longlist}
\end{lemma}
The proof of this result can be found in Section~\ref{secproofs}.

We can strengthen Lemma~\ref{lemmaapprox1} by showing the convergence
above is in fact uniform for $t$ in a compact set.
%
\begin{lemma}
\label{lemmaapprox2}
For $a\in(0,1)$, $\ve>0$ and $M>0$,
\begin{longlist}[(ii)]
\item[(i)]
\[
\lim_{x\to\infty} \mathbb{P} \Bigl(\sup_{t\in[-M,M]}\sup_{u\in[0,a]}
\bigl\llvert Z_0^x\bigl(us_x(t)\bigr)-
\phi_0^x(u,t) \bigr\rrvert >\ve \Bigr) = 0.
\]
\item[(ii)]
\[
\lim_{x\to\infty} \mathbb{P} \Bigl(\sup_{t\in[-M,M]}\sup_{u\in[0,1]}
\bigl\llvert Z_1^x\bigl(us_x(t)\bigr)-
\phi_1^x(u,t) \bigr\rrvert >\ve \Bigr) = 0.
\]
\end{longlist}
\end{lemma}
The result is established by showing that the probabilities in the
statement of Lemma~\ref{lemmaapprox1} are monotone in the parameter
$t$. Again, we defer the full proof until Section~\ref{secproofs}.

Using this uniform approximation result, we establish the following
theorem for the process paths evaluated at multiples of the $Z_0$
extinction time.

\begin{theorem}
\label{thmapprox1}
For $a<1$, $\ve>0$ and $\mu_x = \mu x^{-\alpha}$, where $\alpha\in(0,1)$,
\begin{longlist}[(ii)]
\item[(i)]
\[
\lim_{x\to\infty}\Prob \biggl(\sup_{u\in[0,a]} x^{u-1}\biggl
\llvert Z_0(uT_x)-\phi_0 \biggl(u,T_x-
\frac{1}{r}\log x \biggr)\biggr\rrvert >\ve \biggr) = 0.
\]
\item[(ii)]
\[
\lim_{x\to\infty}\Prob \biggl(\sup_{u\in[0,1]} x^{\alpha-u\lambda
_1/r-1}\biggl
\llvert Z_1(uT_x)-\phi_1
\biggl(u,T_x-\frac{1}{r}\log x \biggr)\biggr\rrvert >\ve \biggr)
= 0.
\]
\end{longlist}
\end{theorem}
\begin{pf}
This result now follows directly from Lemma~\ref{lemmaapprox2} and the
result in \eqref{eqexttimeconv}.
\end{pf}

We define the following stochastic processes: if $u\in[0,1]$,
\begin{eqnarray*}
\psi_0(u)&=&e^{-u(\eta+\log c)},
\\
\psi_1(u)&=&
\cases{\displaystyle\frac{\mu}{\lambda_1+r}e^{({\lambda_1u}/{r})(\eta
+\log c)},&\quad $u>0$,
\vspace*{2pt}
\cr
0,&\quad $u=0$.}
\end{eqnarray*}
These processes represent the limits of our scaled population
processes. However, note that $\psi_1$ is not right-continuous at 0 and
therefore it is\vadjust{\goodbreak} not possible to establish that the scaled population
processes are tight on an interval of the form $[0,b]$ in the standard
Skorokhod topology. This is a result of the massive influx of mutations
near $t=0$ in the unscaled process. 

Next, we utilize Theorem~\ref{thmapprox1} to establish weak
convergence in the Skorokhod sense of the joint process.
%
\begin{corollary}
For $\alpha\in[0,1)$ and $0<a<b<1$ the joint process
\[
\bigl\{\bigl(x^{u-1} Z_0(uT_x),
x^{\alpha- \lambda_1u/r - 1} Z_1(u T_x)\bigr),u\in [a,b]\bigr\}
\Rightarrow\bigl\{ \bigl(\psi_0(u),\psi_1(u) \bigr),u
\in[a,b]\bigr\}
\]
as $x \rightarrow\infty$ in the standard Skorokhod topology, $D([a,b])$.
\end{corollary}
\begin{pf}
For ease of notation, throughout this proof we will use the following notation:
\[
\bphi^x(u)= \biggl(\phi^x_0
\biggl(u,T_x-\frac{1}{r}\log x\biggr),\phi_1^x
\biggl(u,T_x-\frac
{1}{r}\log x\biggr) \biggr).
\]

Clearly, from the result in Theorem~\ref{thmapprox1} it suffices to
prove that as $x\to\infty$
\[
\biggl(\phi_0^x\biggl(\cdot,T_x-
\frac{1}{r}\log x\biggr),\phi_1^x\biggl(
\cdot,T_x-\frac
{1}{r}\log x\biggr) \biggr)\Rightarrow \bigl(
\psi_0(\cdot),\psi(\cdot) \bigr)
\]
in $D([a,b])$. We will carry this out via Theorem 13.3 of \cite
{billingsleyConv}. First, we observe that convergence in finite
dimensional distributions follows from \eqref{eqexttimeconv} and the
continuous mapping theorem. Thus it only remains to establish
tightness. Since our limit functions are continuous at $u=a$ and $u=b$,
it suffices to establish that for every $\ve>0$
%
\begin{equation}
\label{eqtightnesscond} \lim_{\delta\to0}\limsup_{x\to\infty}\Prob\bigl(
\omega_x^{\prime\prime}(\delta)\geq\ve\bigr)=0,
\end{equation}
where
\[
\omega_x^{\prime\prime}(\delta)\doteq\sup\bigl\{\bigl\|
\bphi^x(u)-\bphi^x(u_1)\bigr\|\wedge\bigl\|
\bphi^x(u_2)-\bphi^x(u)\bigr\|\dvtx
u_1\leq u\leq u_2, u_2-u_1
\leq\delta\bigr\}
\]
and $\| x \| \equiv|x_1|+|x_2|$.
From the mean value theorem there exists a constant $C$ such that for
$u<v\in[a,b]$
\begin{eqnarray*}
\bigl\|\bphi^x(u)-\bphi^x(v)\bigr\|&\leq& C(v-u)\biggl|T_x-
\frac{1}{r}\log x\biggr|e^{\lambda_1 v|T_x-({1}/{r})\log x|}
\\
&\leq& C(v-u)e^{(\lambda_1+1)b|T_x-({1}/{r})\log x|}.
\end{eqnarray*}
Thus, if $\omega_x^{\prime\prime}(\delta)\geq\ve$, then
\[
C\delta e^{(\lambda_1+1)b|T_x-({1}/{r})\log x|}\geq\ve,
\]
and therefore,
\[
\Prob\bigl(\omega_x^{\prime\prime}(\delta)\geq\ve\bigr)\leq\Prob
\biggl(\biggl|T_x-\frac
{1}{r}\log x\biggr|\geq b(\lambda+1)\log(\ve/\delta)
\biggr).
\]
Condition \eqref{eqtightnesscond} then follows by taking the
limit as $x\to\infty$ [using \eqref{eqexttimeconv}] and then sending
$\delta$ to 0.
\end{pf}
%
\subsection{Crossover time}
We define the following stochastic time:
\[
\xi\equiv\operatorname{inf} \bigl\{ t >0 | Z_1 (t) \geq
Z_0(t) \bigr\},
\]
which we refer to as the ``crossover'' time, since it is the first time
at which the $Z_0$ and $Z_1$ paths cross, and represents roughly the
time at which the $Z_1$ or resistant cell population begins to dominate
the tumor. In this section we investigate, using the limit theorems
proven in the previous section, the distribution of the crossover time
scaled by $T_x$.
First we utilize the crossover time of the limit processes to obtain an
estimate of this time. In particular, we define $\tilde{u}$ to be the
solution to
\[
\phi_0\biggl(\tilde{u}, T_x - \frac{1}{r} \log x
\biggr) = \phi_1\biggl(\tilde{u}, T_x -
\frac{1}{r} \log x\biggr).
\]
We obtain
%
\[
\tilde{u} = \frac{\log ( \mu+(\lambda_1+r) x^\alpha )-\log\mu
}{T_x(\lambda_1 +r)}
\]
and establish the following result.
%
\begin{theorem}\label{thmcrosstime} The estimate $\tilde{u}$ and the
scaled crossover time, $\xi/T_x$, converge to each other in probability
as $x \rightarrow\infty$, that is, for all $\ve>0,$
\[
\Prob\biggl(\biggl| \frac{\xi}{T_x} - \tilde{u} \biggr|>\ve\biggr)\rightarrow0\qquad \mbox{as
} x \rightarrow\infty.
\]
\end{theorem}
See Section~\ref{secproofs} for the proof of this result, which
follows as an application of Theorem~\ref{thmapprox1}.

\subsection{Turnaround time: Progression of disease}
In this section we characterize the time at which the total tumor
population stops declining and starts increasing. Define the following
set of random times associated with the unscaled escape process:
\[
\tau= \mathop{\operatorname{argmin}}_{t\geq0}\bigl\{Z_0(t) +
Z_1(t) \bigr\}.
\]
Using the sample path approximations, we can approximate this set of
times, rescaled by the extinction time $T_x$, as the random variable
%
\begin{equation}
\label{eqapproxturnaround} u^* \equiv\frac{\log ({r}/({\lambda_1\mu}) )+\log
(x^{\alpha}(\lambda_1+r)-\mu )}{(\lambda_1+r)T_x}.
\end{equation}
This corresponds to the time at which the approximated path of the
total population size has derivative zero.
Looking at the highest order terms in \eqref{eqapproxturnaround}, we
see that for large $x$,
%
\begin{equation}
\label{equ*approx} u^*\approx\frac{\alpha r}{\lambda_1+r}.
\end{equation}
Thus a higher mutation rate, or smaller $\alpha$, leads to a quicker
turnaround time (relative to the extinction time). In addition, as the
decay rate $r$ increases, the time of progression relative to the time
of extinction increases.

Throughout this section we work with the sped-up but unscaled joint
population processes, $Z_i(us_x(t))$.
For simplicity, write the sum of the mean of $Z_0$ and $Z_1$ as
\begin{eqnarray*}
f_{x,t}(u) &\equiv&\EE Z_0\bigl(us_x(t)\bigr) +
\EE Z_1\bigl(us_x(t)\bigr)
\\
&=& xe^{\lambda_0u(({1}/{r})\log x+t)} \biggl(1-\frac{\mu}{x^{\alpha
}(\lambda_1+r)} \biggr)+\frac{x^{1-\alpha}\mu e^{\lambda_1u(
({1}/{r})\log x+t)}}{\lambda_1+r}.
\end{eqnarray*}
We will first show that with high probability, the critical point of $f_{x,t}$,
%
\begin{equation}
\label{eqturnaround} u^\ast(t) \equiv\frac{\log ({r}/{(\lambda_1\mu)} )+\log
(x^{\alpha}(\lambda_1+r)-\mu )}{(\lambda_1+r) (({1}/{r})\log
x+t )}
\end{equation}
is close to the minimum of $Z_0(us_x(t))+Z_1(us_x(t))$. We then
establish that this statement is in fact true uniformly for $t$ in
compact sets, and that $\tau/T_x$ is well approximated by $u^*$.

Since $u^*(t)$ is a critical point of $f_{x,t}$ we have the following
representation that will be useful:
\[
f_{x,t}\bigl(u^\ast(t)\bigr) = xe^{\lambda_0u^*(t)(({1}/{r})\log x+t)} \biggl(1-
\frac{\mu}{x^{\alpha}(\lambda_1+r)} \biggr) \biggl(1+\frac{r}{\lambda_1} \biggr),
\]
and
\begin{eqnarray*}
f_{x,t}\bigl(u^*(t)+y\bigr)
&=& xe^{\lambda_0u^*(t)(({1}/{r})\log x+t)} \biggl(1-\frac{\mu}{x^{\alpha
}(\lambda_1+r)} \biggr)\\
&&{}\times
\biggl(e^{\lambda_0 y(({1}/{r})\log x+t)}+\frac
{r}{\lambda_1}e^{\lambda_1 y(({1}/{r})\log x +t)} \biggr).
\end{eqnarray*}
Therefore,
%
\begin{equation}
\label{eqsteepprop}\qquad f_{x,t}\bigl(u^*(t)+y\bigr) = f_{x,t}
\bigl(u^*(t)\bigr) \biggl[ \biggl(\frac{\lambda_1}{\lambda_1+r} \biggr) \biggl(x^{-y}e^{\lambda_0yt}+
\frac{r}{\lambda_1}x^{\lambda
_1y/r}e^{\lambda_1yt} \biggr) \biggr].
\end{equation}
With this ``steepness'' at the minimum property we can establish that
with high probability (for $x$ large) the minimum of the total
population is achieved at $u^*(t)$.
%
\begin{lemma}
\label{lemturnaroundtime1}
For $\ve> 0$,
\[
\Prob \biggl(\frac{\tau}{({1}/{r})\log x + t}\cap\bigl[u^*(t)-\ve,u^*(t)+\ve
\bigr]^c\neq\varnothing \biggr) \rightarrow0
\]
as $x \rightarrow\infty$.
\end{lemma}
The proof of this result is deferred to Section~\ref{secproofs}.\vadjust{\goodbreak}

Similar to the approximation result in Lemma~\ref{lemmaapprox1}, it is
then possible to establish that an analogous result holds uniformly for
$t$ in compact sets.
%
\begin{lemma}
\label{lemturnaroundunif}
For $\ve> 0$ and a constant $M > 0$,
\begin{longlist}[(ii)]
\item[(i)]
\[
\Prob\Bigl(\sup_{t \in[-M,M]} \inf_{u \in[0, u^*(t) - \ve)}Z_0
\bigl(us_x(t)\bigr) + Z_1\bigl(us_x(t)\bigr)
< Z_0\bigl(u^*T_x\bigr) + Z_1
\bigl(u^*T_x\bigr)\Bigr) \rightarrow0.
\]
\item[(ii)]
\[
\Prob\Bigl(\sup_{t \in[-M,M]} \inf_{u \in[0, u^*(t) + \ve)}Z_0
\bigl(us_x(t)\bigr) + Z_1\bigl(us_x(t)\bigr)
< Z_0\bigl(u^*T_x\bigr) + Z_1
\bigl(u^*T_x\bigr)\Bigr) \rightarrow0
\]
as $x \rightarrow\infty$.
\end{longlist}
\end{lemma}
See Section~\ref{secproofs} for details of the proof.

We can now establish that the turnaround time of the scaled process
$\tau$ normalized by the extinction time $T_x$ converges in probability
to $u^*$.
%
\begin{theorem}
\label{thmturnaroundapprox}
For $\ve> 0$,
\[
\Prob \biggl(\frac{\tau}{T_x}\cap\bigl[u^*-\ve,u^*+\ve\bigr]^c\neq
\varnothing \biggr) \rightarrow0
\]
as $x \rightarrow\infty$.
\label{turnaroundthm}
\end{theorem}
\begin{pf}
Using similar techniques as in the proof of Theorem~\ref{thmapprox1},
the result follows easily from the previous two lemmas.
\end{pf}

\begin{figure}

\includegraphics{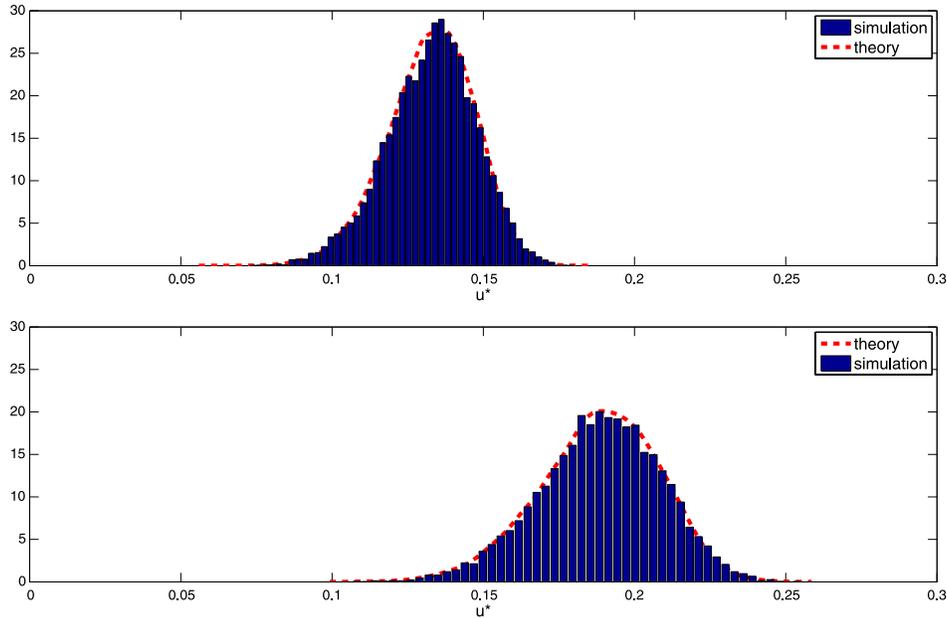}

\caption{Sample PDF of $\tau/T_x$ from simulation of $(Z_0, Z_1)$
process compared with theoretical PDF of $u^*$, for initial size $x =
100\mbox{,}000$ and two parameter sets. Top: $x = 100\mbox{,}000$, $r_0 = 1.0$, $d_0 =
1.5$, $r_1 = 2.0$, $d_1 = 1.0$, $\mu= 0.01$. Bottom: $x = 100\mbox{,}000$, $r_0 =
1.0$, $d_0 = 1.75$, $r_1 = 2.0$, $d_1 = 1.0$, $\mu= 0.01$.}
\label{figturnaroundhisfit}
\end{figure}

In Figure~\ref{figturnaroundhisfit} we compare the sample probability
density function of $\tau/T_x$ from simulations of the $(Z_0, Z_1)$
process with the theoretical PDF of $u^*$. It is observed that even
with an initial starting population of size $x = 100\mbox{,}000$, the
comparisons are favorable. Thus, in the application of interest where
$x$ is on the order of $10^6$ cells or greater, we expect these
limiting approximations to be of use.

\section{An example: Recurrence dynamics in nonsmall cell lung cancer}\label{sec5}

In this section we apply the results to a simple model of drug
resistance in nonsmall cell lung cancer (NSCLC). Nonsmall cell lung
cancer is a disease in which malignant cells form in the tissues of the
lung; it is the most common type of lung cancer, which causes over
150,000 deaths per year in the U.S. In recent years, a new class of
targeted anti-cancer drugs called tyrosine kinase inhibitors has been
developed. These inhibitors target molecules specifically within cancer
cells and inhibit key signaling pathways such as the epidermal growth
factor receptor (EGFR). Two such inhibitors, erlotinib and gefitinib,
have been shown to be extremely successful in reducing tumor burden in
a substantial subset of NSCLC patients. However, point mutations in the
binding site of the drug have been identified that confer resistance to
both therapies, and thus lead to recurrence or progression of the disease.

\begin{figure}

\includegraphics{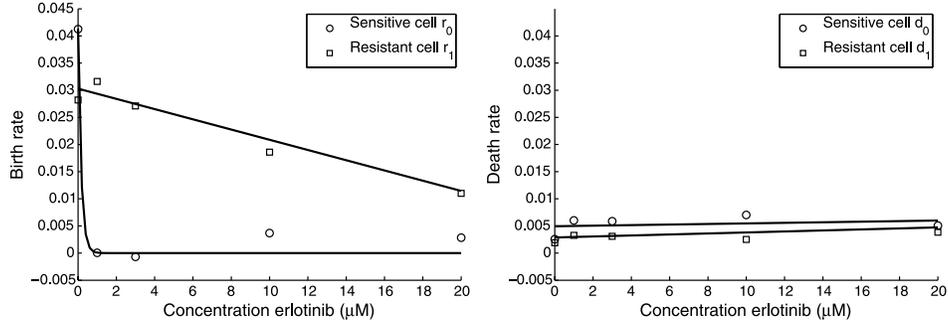}

\caption{Growth and death rate data ($\mbox{hours}^{-1}$) for
erlotinib-sensitive (PC-9) and erlotinib-resistant NSCLC cells as
a function of drug concentration $($data published in~\cite{ChFoSTM2011}).}
\label{figrates}\vspace*{-3pt}
\end{figure}

In previous work~\cite{ChFoSTM2011} we characterized the in vitro
growth rates of a pair of human NSCLC cell lines which were sensitive
or resistant to the drug erlotinib (see Figure~\ref{figrates}).
Here we utilize this experimental growth kinetic data and apply our
results on turnaround time distribution to study the properties of the
time of disease progression. In particular, for a series of drug
concentrations we characterize the distribution of the random time
$u^*$, using the experimental data to ascertain $r_0, d_0, r_1$ and
$d_1$. In addition, we use known estimates of the biological parameter
$\mu_x \approx10^{-8}$, which corresponds to the mutation probability
per cell division per base pair in the genome~\cite{Seshadri87,Oller89}.
We can then apply our estimates of the turnaround time
distribution to study how the time until progression varies as a
function of drug concentration. These distributions of $u^*$ are
helpful in predicting the likely success of the therapy. In particular,
$u^*$ indicates the fraction of the total time that the drug is
effective ($T_x$) at which disease progression occurs. If the
distribution of $u^*$ for a particular drug at a specific concentration
\begin{figure}[b]

\includegraphics{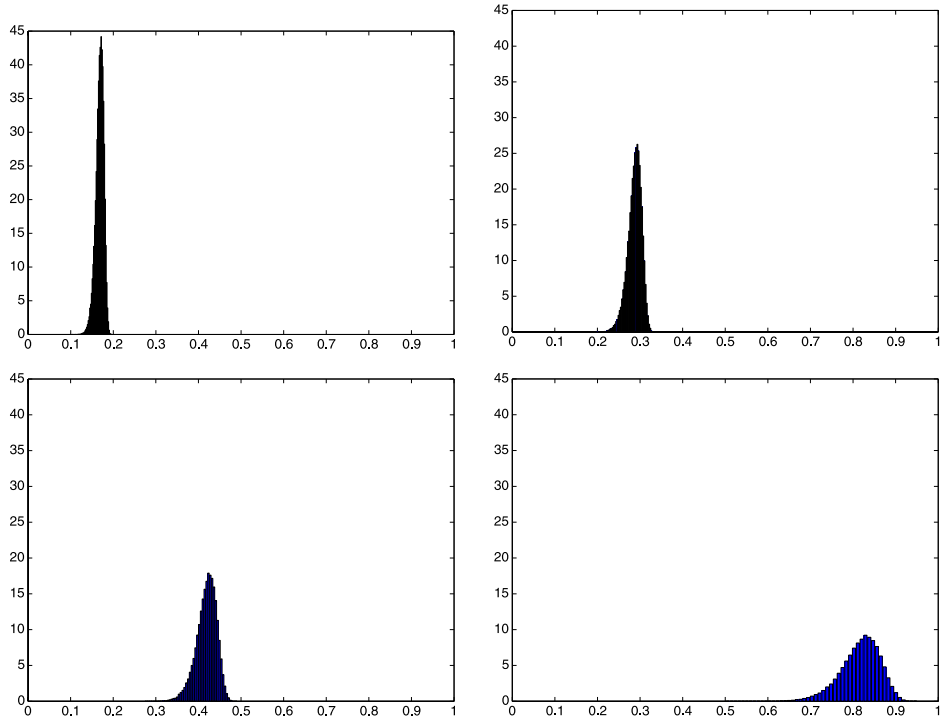}

\caption{Distributions of the turnaround time, $u^*$, for a NSCLC
tumor with initially $10^9$ sensitive cells, treated with erlotinib at
1~$\mu M$ (top left), 3~$\mu M$ (top right), 5~$\mu M$
(bottom left) and 10~$\mu M$ (bottom right).}
\label{figturndist}
\end{figure}
has most of its mass bounded far below 1, the chance that the sensitive
cell drug population is eradicated by the time of progression is
extremely low. On the other hand, drugs whose profiles which place most
of the $u^*$ distribution's mass closer to 1 have better prospects of
eliminating the tumor. In Figure~\ref{figturndist} we plot the $u^*$
distribution for a NSCLC tumor starting with $10^9$ sensitive cells
treated with erlotinib at various concentrations. Note that the current
standard of care, the FDA approved dose elicits a concentration of 3~$\mu M$
in the plasma which corresponds to the upper right plot. As the drug
concentration increases, the distribution of $u^*$ moves accordingly to
the right; however, even at the highest concentration the majority of
the mass is still bounded well below 1 which indicates likely failure
of the therapy. In clinical observations, following an initial response
in terms of tumor reduction, 100 percent of patients develop resistance
usually within 24 months of starting treatment~\cite{ChPa10}.

One major clinical question in NSCLC treatment today is: once the
disease has progressed and the tumor size begins to increase, what
course of therapy is optimal? In particular, should the drug be
withdrawn or should the patient be kept on erlotinib or gefitinib? If
drug is maintained, how long should it be administered beyond
progression? Here, estimates of the $u^*$ distribution can be of use.
We note that $\tau$ is a clinically observable quantity since it
represents the time until disease progression from the start of
treatment. Once $\tau$ is observed, using Theorem~\ref{turnaroundthm}
and the approximation in \eqref{equ*approx} we can approximate $T_x$,
which represents the time at which the entire drug-sensitive population
is eradicated. This gives a clear endpoint, $T_x$ beyond which
erlotinib therapy is unwarranted. Furthermore, we can easily obtain the
distribution of the population size of resistant cells at this time
$Z_1(T_x)$ to estimate the projected resistant tumor size at the time
the sensitive cells are eradicated. This information aids in
determining whether erlotinib treatment should be maintained until
$T_x$ or a switch to alternative therapy should be made prior to
$T_x$.


\section{Summary}\label{sec6}
In this work we have considered the stochastic dynamics of escape from
extinction in a binary branching process model.\vadjust{\goodbreak} By considering the
large starting population limit, we approximate the birth--death process
with a simpler stochastic process whose only randomness is inherited
from the weak limit of the extinction time. Using this limit, we
approximate the distribution of the time until the total population
begins to increase, and the time at which the escape mutants first
begin to dominate the population. One of many possible future
extensions is to consider the problems in this paper in a non-Markovian
setting, that is, nonexponential distribution between events. This work
contributes to a growing body of literature concerned with the
mathematical understanding of cancer evolution, as well as to the
general understanding of extinction and escape paths in branching
process models. In future work we examine the setting $\alpha= 1$,
where $O(1)$ mutations arise before extinction and escape from
extinction is not assured in the large~$x$ limit.

\section{Proof of main results}\label{secproofs}
\subsection{\texorpdfstring{Proof of Theorem \protect\ref{thmmeanscaling}}{Proof of Theorem 2}}
We first establish the scaling of the mean for large initial population
$x$.

By conditioning on the path of $Z_0$ until $T_x$ we get the formula
\eqref{eqexpmuts}; performing a change of measure and
flipping the order of integration (by Tonelli's theorem) we see
\begin{eqnarray*}
\EE\bigl[Z_1(vT_x)\bigr] &=& \mu_x \EE
\biggl[\int_0^{\infty}1_{\{T_x\geq y/v\}}Z_0(y)
\exp \bigl(\lambda_1(vT_x-y) \bigr)\,dy \biggr]
\\
&= &\mu_x\int_0^{\infty}\EE
\bigl[1_{\{T_x\geq y/v\}}Z_0(y)\exp \bigl(\lambda_1(vT_x-y)
\bigr) \bigr]\,dy
\\
&=& \mu_x\int_0^{\infty}\int
_{y/v}^{\infty}\sum_{n=1}^{\infty}ne^{\lambda
_1(vt-y)}
\Prob \bigl(T_x\in dt | Z_0(y)=n \bigr)\Prob
\bigl(Z_0(y)=n\bigr)\,dy.
\end{eqnarray*}
Note that $\Prob(T_x\in dt | Z_0(y)=n) = g_n(t-y)\,dt$, where $g_n$ is
the density of the extinction time for a population starting from a
population of size $n$, and can be written as $g_n(t) = n
(G(t) )^{n-1}g(t)$, where $g$ is the density of the extinction
time for a population starting from a single cell, and $G$ is the c.d.f.
Therefore, upon rearranging the order of integration we get that
\begin{eqnarray*}
&&\hspace*{-3pt}\EE\bigl[Z_1(vT_x)\bigr]
\\
&&\hspace*{-3pt}\qquad = \mu_x\int_0^{\infty}\int
_{y/v}^{\infty}e^{\lambda
_1(vt-y)}g(t-y) \Biggl(\sum
_{n=1}^{\infty}n^2 \bigl(G(t-y)
\bigr)^{n-1}\Prob\bigl(Z_0(y)=n\bigr) \Biggr)\,dt\,dy.
\end{eqnarray*}
Next define
\[
F_x(s,t)=\EE\bigl[ s^{Z_0(t)}\bigr]\quad\mbox{and}\quad F(s,t)=\EE
\bigl[s^{Z_0(t)}|Z_0(0)=1,Z_1(0)=0\bigr],
\]
and observe that due to the independence of the branching structure
$F_x(s,t) =  (F(s,t) )^x$.
Therefore,
%
\begin{eqnarray}
\label{eqgenfun}
&&\sum_{n=1}^{\infty} n^2s^{n-1}
\Prob\bigl(Z_0(t)=n\bigr)
\nonumber\\
&&\qquad=s\frac{\partial^2}{\partial s^2}F_x(s,t)+\frac{\partial}{\partial
s}F_x(s,t)
\nonumber
\\[-8pt]
\\[-8pt]
\nonumber
&&\qquad=sx \bigl(F(s,t) \bigr)^{x-1}\frac{\partial^2}{\partial
s^2}F(s,t)+sx(x-1)
\bigl(F(s,t) \bigr)^{x-2} \biggl(\frac{\partial
}{\partial s}F(s,t)
\biggr)^2
\nonumber
\\
&&\qquad\quad{} +x \bigl(F(s,t) \bigr)^{x-1}\frac{\partial}{\partial s}F(s,t).\textbf{\nonumber}
\end{eqnarray}
For ease of notation we will simply write $\frac{\partial}{\partial
s}F(s,t) = F'(s,t)$.
Using \eqref{eqgenfun} we obtain
\begin{eqnarray*}
&&\EE\bigl[Z_1(vT_x)\bigr]
\\
&&\qquad=x\mu_x\int_0^{\infty}\int
_{{y}/{v}}^{\infty}e^{\lambda
_1(vt-y)}g(t-y) F\bigl(G(t-y),y
\bigr)^{x-2}\\
&&\hspace*{92pt}{}\times\bigl[ (x-1)G(t-y)F'\bigl(G(t-y),y
\bigr)^2
\\
&&\hspace*{88pt}\qquad{} + G(t-y)F\bigl(G(t-y),y\bigr)F''\bigl(G(t-y),y
\bigr)\\
&&\hspace*{130pt}\qquad{}+F\bigl(G(t-y),y\bigr)F'\bigl(G(t-y),y\bigr)\bigr]\,dy\,dt.
\end{eqnarray*}
This expression can be analyzed using techniques from \cite
{JaKlSa07Adv}. In particular, if we introduce the change of variable
\[
t = \frac{1}{r} (z+\log cx ) = z_x,
\]
observe that
\[
e^{\lambda_1(vz_x-y)} = e^{\lambda_1(vz/r-y)}(cx)^{\lambda_1 v/r}.
\]
After the change of variables we have
\begin{eqnarray*}
\EE\bigl[Z_1(vT_x)\bigr] =I_1(x,v)+I_2(x,v),
\nonumber
\end{eqnarray*}
where
\begin{eqnarray*}
I_1(x,v) &=& \frac{x(x-1)(cx)^{{\lambda_1v}/{r}}\mu_x}{r}\int_0^{\infty}
\int_{{ry}/{v}-\log cx} \Phi_1(z,y) \,dz\,dy,
\\
I_2(x,v)&= &\frac{x(cx)^{{\lambda_1v}/{r}}\mu_x}{r}\int_0^{\infty
}
\int_{{ry}/{v}-\log cx}^{\infty} \Phi_2(z,y) \,dz\,dy
\end{eqnarray*}
and
\begin{eqnarray*}
\Phi_1(z,y)
&=& e^{\lambda_1(vz/r-y)}g(z_x-y)G(z_x-y)\\
&&{}\times \bigl(F
\bigl(G(z_x-y),y\bigr) \bigr)^{x-2} \bigl(F'
\bigl(G(z_x-y),y\bigr) \bigr)^2,
\\
\Phi_2(z,y)
&=&e^{\lambda_1(vz/r-y)}g(z_x-y) \bigl(F\bigl(G(z_x-y),y\bigr)
\bigr)^{x-1}\\
&&{}\times \bigl[G(z_x-y)F''
\bigl(G(z_x-y),y\bigr) + F'\bigl(G(z_x-y),y\bigr)\bigr].
\end{eqnarray*}

We will now establish that for $v\in(0,1]$, $x^{\alpha-1-{\lambda_1v}/{r}}I_1(x,v)=\tilde{I}_1(x,v)\to I_1(v)$ and $x^{\alpha-1-
{\lambda_1v}/{r}}I_2(x,v)=\tilde{I}_2(x,v)\to0$ as $x\to\infty$. The
integrand of $\tilde{I}_1(x,v)$ is
\begin{eqnarray*}
f_x(z,y)
&=& c_0xg(z_x-y)G(z_x-y) \\
&&{}\times\bigl(F
\bigl(G(z_x-y),y\bigr) \bigr)^{x-2} \bigl(F'
\bigl(G(z_x-y),y\bigr) \bigr)^2 e^{\lambda_1(vz/r-y)},
\end{eqnarray*}
where $c_0 = \mu c^{\lambda_1v/r}/r$.
From~\cite{JaKlSa07Adv} we know that as $z\to\infty$, $g(z)\sim
rce^{-rz}$, and therefore
%
\begin{equation}
\label{eqgasymp} g(z_x-y)\sim\frac{r}{x}e^{ry}e^{-z}
\end{equation}
as $x\to\infty$. Next note that there exists a $\xi_x\in(G(z_x-y),1)$
such that
\begin{eqnarray*}
F'\bigl(G(z_x-y),y\bigr) &=& F'(1,y) +
\bigl(1-G(z_x-y) \bigr)F''(
\xi_x,y)
\\
& =& e^{-ry} +O \bigl(1-G(z_x-y) \bigr),
\end{eqnarray*}
and therefore
%
\begin{equation}
\label{eqFprimeasymp} F'\bigl(G(z_x-y),y\bigr)\sim
e^{-ry}.
\end{equation}
Last, observe that
\[
F\bigl(G(z_x-y),y\bigr) = 1 +e^{-ry} \bigl(G(z_x-y)-1
\bigr)+ O \bigl(G(z_x-y)-1 \bigr)^2,
\]
and therefore
\[
\log F\bigl(G(z_x-y),y\bigr) = - \bigl(1-G(z_x-y)
\bigr)e^{-ry}+O \bigl(G(z_x-y)-1 \bigr)^2.
\]
Observe that
\[
1-G(z_x-y)\sim\frac{e^{ry}e^{-z}}{x},
\]
which gives that
%
\begin{equation}
\label{eqFpowerasymp} \bigl(F\bigl(G(z_x-y),y\bigr)
\bigr)^{x-2}\sim\exp\bigl[-e^{-z}\bigr].
\end{equation}
Combining \eqref{eqgasymp}, \eqref{eqFprimeasymp} and \eqref
{eqFpowerasymp} we see that
%
\begin{equation}
\label{eqintegrandlimit} \lim_{x\to\infty}f_x(z,y) =
c_1 e^{-z(1-({\lambda
_1}/{r})v)}e^{-y(r+\lambda_1)}\exp\bigl[-e^{-z}
\bigr],
\end{equation}
where $c_1=rc_0$.
In order to evaluate the limit of $\tilde{I}_1$ it thus remains to show
that the limit can be passed inside the integral; this will be done by
finding an integrable function $h$ such that $f_x(z,y)\leq h(z,y)$.
First note that since $G(z)\leq1$ and $F'(s,t)\leq E_1 Z_0(t) =
e^{-rt}$, we have
\[
f_x(z,y)\leq c_0x e^{\lambda_1(vz/r-y)}g(z_x-y)
\bigl(F\bigl(G(z_x-y),y\bigr) \bigr)^{x-2}e^{-2ry}.
\]
Then observe that there exists a constant $k_1$ such that
%
\begin{equation}
\label{eqgupper} g(z_x-y)\leq k_1e^{-r(z_x-y)} =
\frac{k_2}{x}e^{-z}e^{ry}.
\end{equation}
Since $\log x\leq x-1$ we have
\begin{eqnarray*}
\bigl(F\bigl(G(z_x-y),y\bigr) \bigr)^{x-2} &=& \exp
\bigl[(x-2)\log F\bigl(G(z_x-y),y\bigr) \bigr]
\\
& \leq&\exp \bigl[-(x-2) \bigl(1-F\bigl(G(z_x-y),y\bigr)\bigr)
\bigr].
\end{eqnarray*}
Using results from the proofs of Proposition 1 and Lemma 1 of \cite
{JaKlSa07Adv}, we can establish that
%
\begin{equation}
\label{eqFpowerupper} 1 - F\bigl(G(z_x-y),y\bigr) \geq
\frac{e^{-z}}{x}.
\end{equation}
Based on \eqref{eqgupper} and \eqref{eqFpowerupper} we see that we
can use the dominating function
\[
h(z,y)=k_3e^{\lambda_1({vz}/{r}-y)}e^{-z}e^{-ry}\exp
\bigl[-k_2e^{-z}\bigr].
\]
With this result we see that
\begin{eqnarray*}
\lim_{x\to\infty}\tilde{I}_1(x,v) &=& c_0\int
_0^{\infty}\int_{-\infty
}^{\infty}
\lim_{x\to\infty}f_x(z,y)\,dz\,dy
\\
&=& c_1\int_0^{\infty}\int
_{-\infty}^{\infty}e^{\lambda_1(
{vz}/{r}-y)}e^{-yr}e^{-z}
\exp\bigl[-e^{-z}\bigr]\,dz\,dy
\\
&=& \frac{\mu c^{\lambda_1v/r}\Gamma(1-{\lambda_1v}/{r})}{(\lambda_1+r)}.
\end{eqnarray*}
We now consider $\tilde{I}_2$. First observe that for $(s,t)\in
[0,1]\times[0,\infty)$ there exists finite~$k_4$ such that
$F'(s,t)\leq k_4$ and $F''(s,t)\leq k_4$ and, of course, $F(s,t)\leq
1$. Therefore, if we consider $\tilde{I}_2$ in terms of the original
variables, there exists a finite constant $k_5$ such that
\begin{eqnarray*}
\tilde{I}_2(x,v) &\leq& k_5x^{-v\lambda_1/r}\int
_0^{\infty}\int_{y/v}^{\infty}e^{\lambda_1(vt-y)}g(t-y)\,dt
\,dy
\\
&=& k_5x^{-v\lambda_1/r}\int_0^{\infty}
\int_{y(1/v-1)}^{\infty}e^{-\lambda
_1y}g(s)e^{\lambda_1v(s+y)}\,ds\,dy
= k_6x^{-v\lambda_1/r},
\end{eqnarray*}
where the first equality follows by using the change of variable
$s=t-y$. Thus $\tilde{I}_2(x,v)\to0$ as $x\to\infty$ for $v\in(0,1]$.

\subsection{\texorpdfstring{Proof of Lemma \protect\ref{lemmamoments}}{Proof of Lemma 1}}

We will start by establishing item (ii). Define $\mathcal{F}_{\infty
}^0$ to be the sigma algebra generated by the wave 0 population until
their eventual extinction, then
\[
E \bigl[Z_1\bigl(s_x(t)\bigr) |\mathcal{F}_{\infty}^0
\bigr] = \mu_x\int_0^{s_x(t)}Z_0(y)e^{\lambda_1(s_x(t) -y)}\,dy
\]
and therefore
\[
\EE \bigl[Z_0(s)Z_1\bigl(s_x(t)\bigr)
\bigr] = e^{\lambda_1 s_x(t)}\mu_x\int_0^{s_x(t)}
\EE\bigl[Z_0(y)Z_0(s)\bigr]e^{-\lambda_1 y}\,dy.
\]

Now we establish the second moment result, item (i). For simplicity we
evaluate $\EE[Z_1(t)]$ for a positive $t$. For ease of notation we will
use the following $\tilde{\EE}[\cdot] = \EE[\cdot|\mathcal{F}_{\infty
}^0]$. Consider a partition of $[0,t]$, $0<\Delta<2\Delta<\cdots<t$,
where $\Delta=t/m$ for a large integer $m$. Then we can write
\begin{eqnarray*}
\EE\bigl[Z_1(t)^2|\mathcal{F}_{\infty}^0
\bigr] &=& \tilde\EE \Biggl(\sum_{j=0}^{m}
\sum_{k=1}^{N_j}B_{j,k}(t-
\tau_{j,k}) \Biggr)^2
\\
&=& \sum_{j=0}^m \tilde\EE \Biggl(\sum
_{k=1}^{N_j} B_{j,k}(t-
\tau_{j,k}) \Biggr)^2
\\
&&{} + \sum_{j=0}^m\sum
_{\ell\neq j}\tilde\EE \Biggl(\sum_{k=1}^{N_j}B_{j,k}(t-
\tau_{j,k}) \Biggr)\tilde\EE \Biggl(\sum_{k=1}^{N_{\ell}}B_{\ell,k}(t-
\tau_{\ell,k}) \Biggr),
\end{eqnarray*}
where $N_j$ is the number of type-1 mutants created in $[j\Delta
,(j+1)\Delta)$, $\{B_{j,k}\}$ is a collection of i.i.d binary
birth--death processes with birth rate $a_1$ and death rate $d_1$ and
$\tau_{j,k}$ is the time of creation for the $k$th mutant created in
$[j\Delta,(j+1)\Delta)$. In the previous display we have used the
independence of the branching process to derive the second equality.
For $0\leq j\leq m$,
\begin{eqnarray*}
\EE N_j &=& \Delta\mu_x Z_0(j\Delta)+o(
\Delta),
\\
\EE N_j^2 &=&\Delta\mu_x Z_0(j
\Delta) \bigl(1 + \Delta\mu_x Z_0(j\Delta )\bigr)+o(
\Delta).
\end{eqnarray*}
Therefore,
\begin{eqnarray*}
\tilde\EE \Biggl(\sum_{k=1}^{N_j}
B_{j,k}(t-\tau_{j,k}) \Biggr)^2 &=& \tilde \EE
N_j \tilde\EE B(t-\tau_j)^2 + \tilde\EE
\bigl[N_j(N_j-1)\bigr] \bigl(\tilde\EE B(t-
\tau_j)\bigr)^2
\\
&=& \Delta\mu_x Z_0(j\Delta)\tilde\EE B(t-
\tau_j)^2 + O\bigl(\Delta^2\bigr)
\\
&=& \Delta\mu_x Z_0(j\Delta) \biggl(\frac{2r_1}{\lambda_1}e^{2\lambda
_1(t-\tau_j)}-
\frac{r_1+d_1}{\lambda_1}e^{\lambda_1(t-\tau_j)} \biggr)
\nonumber
\\
&&{} + O\bigl(\Delta^2\bigr)
\end{eqnarray*}
and
\begin{eqnarray*}
\tilde\EE \Biggl(\sum_{k=1}^{N_{\ell}}B_{\ell,k}(t-
\tau_{\ell,k}) \Biggr)&=&\Delta\mu_x Z_0(j\Delta)
\tilde\EE B(t-\tau_j)
\\
&= &\Delta\mu_x Z_0(j\Delta)e^{\lambda_1(t-\tau_j)}.
\end{eqnarray*}
Using the previous two expressions we get
\begin{eqnarray*}
\tilde\EE\bigl[Z_1(t)^2\bigr] &=& \Delta
\mu_x\sum_{j=0}^m
Z_0(j\Delta) \biggl(\frac
{2r_1}{\lambda_1}e^{2\lambda_1(t-\tau_j)}-
\frac{r_1+d_1}{\lambda_1}e^{\lambda_1(t-\tau_j)} \biggr)
\\
&&{} + (\mu_x\Delta)^2\sum_{j=0}^m
\sum_{\ell=0,\ell\neq j}^mZ_0(j\Delta
)Z_0(\ell\Delta)e^{\lambda_1(t-\tau_{j})}e^{\lambda_1(t-\tau_{\ell})}.
\end{eqnarray*}
Sending $\Delta\to0$, integrating over $Z_0$ and replacing $t$ with
$\frac{1}{r}\log x +t$ gives us the desired formula for item (i).

Item (iii) follows immediately from item (i).

\subsection{\texorpdfstring{Proof of Lemma \protect\ref{lemmaapprox1}}{Proof of Lemma 2}}
We will prove the more difficult second statement first. We observe
that it suffices to prove that as $x\to\infty$,
\[
\Prob \Bigl(\sup_{u\in[0,1]} e^{-\lambda_1ut}\bigl\llvert
Z_1^x\bigl(us_x(t)\bigr)-
\phi_1^x(u,t) \bigr\rrvert >e^{-\lambda_1t}\ve \Bigr)
\rightarrow0.
\]
Next observe that
\begin{eqnarray*}
&&\frac{e^{-\lambda_1 ut}}{x^{-\alpha+ 1+\lambda_1u/r}} \bigl( Z_1\bigl(us_x(t)\bigr)-
\phi_1(u,t) \bigr)
\\
&&\qquad= \biggl(\frac{e^{-\lambda_1 ut}}{x^{ 1+\lambda_1u/r}}Z_1\bigl(us_x(t)\bigr)-
\frac
{\mu}{x^{1+\alpha}}\int_0^{u s_x(t)}Z_0(s)e^{-\lambda_1 s}\,ds
\biggr) x^{\alpha}
\\
&&\qquad\quad{} + \frac{\mu}{x}\int_0^{u s_x(t)}
\bigl(Z_0(s)-xe^{\lambda_0 s}\bigr)e^{-\lambda_1s}\,ds,
\end{eqnarray*}
and therefore
\begin{eqnarray*}
&& \Prob \Bigl(\sup_{u\in[0,1]} e^{-\lambda_1ut}\bigl\llvert
Z_1^x\bigl(us_x(t)\bigr)-
\phi_1^x(u,t) \bigr\rrvert >e^{-\lambda_1t}\ve \Bigr)
\\
&&\qquad\leq \Prob \biggl(\sup_{u\in[0,1]} x^{\alpha} \biggl\llvert
\frac{e^{-\lambda_1
ut}}{x^{ 1+\lambda_1u/r}}Z_1\bigl(us_x(t)\bigr)-
\frac{\mu}{x^{1+\alpha}}\int_0^{u
s_x(t)}Z_0(s)e^{-\lambda_1 s}\,ds
\biggr\rrvert >\ve/2 \biggr)
\\
&&\qquad\quad{} + \Prob \biggl(\sup_{u\in[0,1]}\frac{\mu}{x}\int_0^{u
s_x(t)}\bigl|Z_0(s)-xe^{\lambda_0 s}\bigr|e^{-\lambda_1s}\,ds>
\ve/2 \biggr).
\end{eqnarray*}
However, we can observe that the process considered in the second
expression in the sum is monotonic in $u$, and the process considered
in the first expression is a martingale in $u$, which allows us to
arrive at the following simpler inequality:
%
\begin{eqnarray}
\label{eqapprox1ineq1}
&& \Prob \Bigl(\sup_{u\in[0,1]} e^{-\lambda_1ut}\bigl
\llvert Z_1^x\bigl(us_x(t)\bigr)-
\phi_1^x(u,t) \bigr\rrvert >e^{-\lambda_1t}\ve \Bigr)
\nonumber\\
&&\qquad\leq \frac{4x^{2\alpha}}{\ve^2} \EE \biggl[ \biggl(\frac{e^{-\lambda
_1t}}{x^{\lambda_1/r+1}}Z_1
\bigl(s_x(t)\bigr)-\frac{\mu}{x^{1+\alpha}}\int_0^{s_x(t)}Z_0(s)e^{-\lambda_1s}\,ds
\biggr)^2 \biggr]
\\
&&\qquad\quad{} + \Prob \biggl(\frac{\mu}{x}\int_0^{s_x(t)}\bigl|Z_0(s)-xe^{\lambda_0
s}\bigr|e^{-\lambda_1s}\,ds>
\ve/2 \biggr).\nonumber
\end{eqnarray}
Consider the latter quantity first, where it suffices to show that as
$x\rightarrow\infty$,
\[
\frac{\mu}{x}\int_0^{s_x(t)}\EE \bigl[\bigl
\llvert Z_0(s)-xe^{\lambda_0
s}\bigr\rrvert \bigr]e^{-\lambda_1 s}\,ds
\rightarrow0.
\]
%
Next observe that
\[
\operatorname{Var}\bigl(Z_0(s)\bigr) = x \biggl(\frac{r_0+d_0}{\lambda_0}
\biggr) \bigl(e^{2\lambda_0 s}-e^{\lambda_0 s} \bigr).
\]
It follows from the Cauchy--Schwarz inequality that
\[
\frac{\mu}{x}\int_0^{s_x(t)}\EE \bigl[\bigl
\llvert Z_0(s)-xe^{\lambda_0
s}\bigr\rrvert \bigr]e^{-\lambda_1 s}\,ds
= O\bigl(x^{-1/2}\bigr).
\]
Moving on to the first term in \eqref{eqapprox1ineq1},
%
\begin{eqnarray}
\label{eqapprox1eq1}
&& x^{2\alpha} \EE \biggl[ \biggl(\frac{e^{-\lambda_1t}}{x^{\lambda
_1/r+1}}Z_1
\bigl(s_x(t)\bigr)-\frac{\mu}{x^{1 + \alpha}}\int_0^{s_x(t)}Z_0(s)e^{-\lambda_1s}\,ds
\biggr)^2 \biggr]\nonumber
\\
&&\qquad= \biggl(\frac{e^{-\lambda_1 t}x^{\alpha}}{x^{1+\lambda_1/r}} \biggr)^2\EE \bigl[Z_1
\bigl(s_x(t)\bigr)^2 \bigr]
\nonumber
\\[-8pt]
\\[-8pt]
\nonumber
&&\qquad\quad{} - 2\frac{x^{2\alpha}\mu e^{-\lambda_1 t}}{x^{2+\lambda_1 /r + \alpha
}}\int
_0^{s_x(t)}e^{-\lambda_1 s}\EE \bigl[Z_0(s)Z_1
\bigl(s_x(t)\bigr) \bigr]\,ds
\\
&&\qquad\quad{} + \biggl(\frac{\mu}{x} \biggr)^2\int_0^{s_x(t)}
\int_0^{s_x(t)}\EE \bigl[Z_0(s)Z_0(y)
\bigr]e^{-\lambda_1s}e^{-\lambda_1 y}\,ds \,dy.
\nonumber
\end{eqnarray}
Using Lemma~\ref{lemmamoments} we see that \eqref{eqapprox1eq1} can
be written as
\begin{eqnarray*}
\frac{2 r_1 x^{\alpha-1}\mu}{\lambda_1} \int_0^{s_x(t)}
e^{s(\lambda
_0-2\lambda_1)} \,ds - \frac{(r_1 + d_1)\mu e^{-\lambda_1 t}}{\lambda_1
x^{1-\alpha+ \lambda_1/r}} \int_{0}^{s_x(t)}e^{s(\lambda_0-\lambda_1)}\,ds
= O \bigl(x^{\alpha-1} \bigr),
\end{eqnarray*}
thus establishing the result (ii).

We now move on to the proof of item (i). First observe that
\begin{eqnarray*}
\sup_{u\in[0,a]}\bigl\llvert Z_0^x
\bigl(us_x(t)\bigr)-\phi_0^x(u,t)\bigr\rrvert
\leq e^{-\lambda_0at^-}\sup_{u\in[0,a]}e^{-\lambda_0ut}\bigl\llvert
Z_0^x\bigl(us_x(t)\bigr)-
\phi_0^x(u,t)\bigr\rrvert ,
\end{eqnarray*}
where $t^-=-\min(t,0)$.
We will show that as $x\to\infty$,
%
\[
\Prob \Bigl(\sup_{u\in[0,a]}\bigl\llvert e^{-\lambda_0ut}
Z_0^x\bigl(us_x(t)\bigr)-1\bigr\rrvert \geq
\ve e^{\lambda_0 a t^-} \Bigr) \rightarrow0. 
\]
Observe that $e^{-\lambda_0ut} x^{u-1}Z_0 (us_x(t)) -1$ is a
martingale with respect to $u$. Therefore, it suffices to show that as
$x\to\infty$,
\begin{eqnarray*}
&\EE \bigl[e^{-\lambda_0at} Z_0^x\bigl(as_x(t)
\bigr)-1 \bigr]^2 
= x^{2a-1}e^{-2a\lambda_0t}\operatorname{Var}_1Z_0
\bigl(as_x(t)\bigr)\to0,
\end{eqnarray*}
where $\operatorname{Var}_1 Z_0(t)$ represents the variance of $Z_0(t)$
starting with an initial population size 1. The previous expression
reduces to
\[
x^{a-1}e^{-a\lambda_0t}\bigl(x^{-a}e^{-\lambda_0at}-1
\bigr) (r_0+d_0)/\lambda_0,
\]
and since we have assumed that $a<1$, the result is established.

\subsection{\texorpdfstring{Proof of Lemma \protect\ref{lemmaapprox2}}{Proof of Lemma 3}}
Throughout the proof assume that $x>e^{rM}$.
We first establish the result for the $Z_1$ population by showing the
following monotonicity property in $t$:
%
\begin{eqnarray}
\label{eqlemmaapprox2mono1}
&&\sup_{u\in[0,1]}\biggl\llvert
Z_1^x \biggl(u \biggl(\frac{1}{r} \log x +
t_0 \biggr) \biggr)-\phi_1^x
(u,t_0 )\biggr\rrvert
\nonumber
\\[-8pt]
\\[-8pt]
\nonumber
&&\qquad\leq \sup_{u\in[0,1]}\biggl\llvert Z_1^x \biggl(u
\biggl(\frac{1}{r} \log x + t_1 \biggr) \biggr)-
\phi_1^x (u,t_1 ) \biggr\rrvert
\end{eqnarray}
for $t_0\leq t_1$. For any $u\in[0,1]$ set
\[
\bar{u} \equiv\frac{u (({1}/{r})\log x+t_0 )}{
({1}/{r})\log x+t_1},
\]
which, of course, implies that $\bar u (\frac{1}{r}\log x+t_1
) = u (\frac{1}{r}\log x+t_0 )$.
In addition, observe that $\bar u\leq u$ and thus $x^{-\lambda
_1u/r}\leq x^{-\lambda_1\bar u/r}$.
Therefore,
\begin{eqnarray*}
&&x^{-\lambda_1u/r}\biggl\llvert Z_1 \biggl(u \biggl(
\frac{1}{r} \log x + t_0 \biggr) \biggr)-\phi_1
(u,t_0 )\biggr\rrvert
\\
&&\qquad\leq x^{-\lambda_1\bar u/r}\biggl\llvert Z_1 \biggl(\bar u \biggl(
\frac{1}{r}\log x+t_1 \biggr) \biggr)-\phi_1 (
\bar u,t_1 )\biggr\rrvert .
\end{eqnarray*}
Since $u\in[0,1]$ was arbitrary we have that
\begin{eqnarray*}
&&\sup_{t\in[-M,M]}\sup_{u\in[0,1]}x^{\alpha-1-\lambda_1u/r}\bigl\llvert
Z_1 \bigl(us_x(t) \bigr)-\phi_1 (u,t )
\bigr\rrvert
\\
&&\qquad\leq \sup_{u\in[0,1]}x^{\alpha-1-\lambda_1u/r}\biggl\llvert Z_1
\biggl(u \biggl(\frac
{1}{r}\log x+M \biggr) \biggr)-\phi_1
(u,M )\biggr\rrvert .
\end{eqnarray*}
Result (ii) now follows by an application of Lemma~\ref{lemmaapprox1}.

The proof of result (i) will follow a similar approach. In particular,
for $t\leq M$ and $u\in[0,a]$, define
\[
\hat u = u \biggl(\frac{({1}/{r})\log x + t}{({1}/{r})\log x +
M} \biggr).
\]
Notice that
\[
u-\hat u = u \biggl(\frac{M-t}{({1}/{r})\log x + M} \biggr)\leq a \biggl(\frac{M-t}{({1}/{r})\log x + M}
\biggr) = n(x,M).
\]
Using the definition of $n(x,M)$ it follows that $n(x,M)\log x\leq
2arM$ which implies that $x^{u-\hat u}\leq e^{2arM}$.
Based on the definition of $\hat u$ and the upper bound on $x^{u-\hat u}$
\begin{eqnarray*}
&&x^u\bigl\llvert Z_0\bigl(us_x(t)\bigr)-
\phi_0(u,t)\bigr\rrvert
\\
&&\qquad= x^{u-\hat u}x^{\hat u}\biggl\llvert Z_0^x
\biggl(\hat u \biggl(\frac{1}{r}\log x+M \biggr) \biggr)-
\phi_0(\tilde u,M)\biggr\rrvert
\\
&&\qquad\leq e^{2arM}x^{\hat u}\biggl\llvert Z_0 \biggl(
\hat u \biggl(\frac{1}{r}\log x+M \biggr) \biggr)-\phi_0(\hat
u,M)\biggr\rrvert .
\end{eqnarray*}
Since the previous inequality holds for any $u$, we know that for any
$t\in[-M,M]$,
\begin{eqnarray*}
&&\sup_{u\in[0,a]}x^u\bigl\llvert Z_0
\bigl(us_x(t)\bigr)-\phi_0(u,t)\bigr\rrvert
\\
&&\qquad\leq e^{2arM}\sup_{u\in[0,a]}x^{u}\biggl\llvert
Z_0 \biggl(u \biggl(\frac{1}{r}\log x+M \biggr) \biggr)-
\phi_0(u,M)\biggr\rrvert .
\end{eqnarray*}
Thus the result of (i) is established by using the result of Lemma \ref
{lemmaapprox1} for $t=M$.

\subsection{\texorpdfstring{Proof of Theorem \protect\ref{thmcrosstime}}{Proof of Theorem 5}}

First we prove that
$
\Prob(\xi/T_x \leq\tilde{u} - \epsilon)\rightarrow0
$
as $x \rightarrow\infty$. In particular, recall that
\begin{eqnarray*}
\phi_0(u,t) &= &x^{1-u}e^{\lambda_0 u t},
\\
\phi_1(u,t) &=&\frac{\mu x^{1-\alpha+\lambda_1 u/r}e^{\lambda
_1ut}}{\lambda_1+r} \bigl(1-e^{u(\lambda_0-\lambda_1)t}x^{(\lambda
_0-\lambda_1)u/r}
\bigr).
\end{eqnarray*}
Then let us utilize the notation $d(T_x) \equiv T_x - \frac{1}{r} \log
x$ to represent the deviation of~$T_x$ from its scaling,
\begin{eqnarray*}
&&\hspace*{-3pt}\Prob \Bigl( \sup_{u \leq\tilde{u} - \epsilon} \bigl(Z_1(uT_x) -
Z_0(uT_x) \bigr) > 0 \Bigr)
\\
&&\hspace*{-3pt}\qquad \leq\Prob \biggl( \sup_{u \leq\tilde{u} - \epsilon} x^{u-1}(Z_1(uT_x)
- \phi_1\bigl(u, d(T_x)\bigr) + \bigl(
\phi_1\bigl(u, d(T_x)\bigr) - \phi_0\bigl(u,
d(T_x)\bigr)\bigr)
\\
&&\hspace*{75pt}\qquad\quad\phantom{\Prob(}{} + \bigl(\phi_0 \bigl( u, d(T_x)\bigr) -
Z_0^x(uT_x) \bigr) > 0, \bigl|d(T_x)\bigr|
\leq\frac{1}{r}\log x \biggr)
\\
&&\hspace*{-3pt}\qquad\quad{} +\Prob \biggl(\bigl|d(T_x)\bigr|>\frac{1}{r}\log x \biggr)
\\
&&\hspace*{-3pt}\qquad\leq\Prob\biggl( \sup_{u \leq\tilde{u} - \epsilon} x^{u-1}\bigl(Z_1(uT_x)
- \phi_1\bigl(u, d(T_x)\bigr)\bigr)
\\
&&\hspace*{-3pt}\qquad\quad\phantom{\Prob(}{} + \sup_{u \leq\tilde{u} - \epsilon} x^{u-1} \bigl(\phi_1\bigl(u,
d(T_x)\bigr) - \phi_0\bigl(u, d(T_x)\bigr)
\bigr)
\\
&&\hspace*{-3pt}\qquad\quad\phantom{\Prob(}{} + \sup_{u \leq\tilde{u} - \epsilon} x^{u-1}\bigl(\phi_0 \bigl( u,
d(T_x)\bigr) - Z_0^x(uT_x) \bigr)
> 0, \bigl|d(T_x)\bigr|\leq\frac{1}{r}\log x \biggr)
\\
&&\hspace*{-3pt}\qquad\quad{} +\Prob \biggl(\bigl|d(T_x)\bigl|>\frac{1}{r}\log x \biggr).
\end{eqnarray*}
Clearly $\Prob (|d(T_x)|>\frac{1}{r}\log x ) \to0$ as $x\to
\infty$, and it thus remains to analyze the first expression on the
right-hand side of previous display. First notice that if $|d(T_x)|\leq
\frac{1}{r}\log x$, then
\begin{eqnarray*}
&&\sup_{u \leq\tilde{u} - \epsilon} x^{u-1} \bigl(\phi_1\bigl(u,
d(T_x)\bigr) - \phi_0\bigl(u, d(T_x)\bigr)
\bigr)
\\
&&\qquad= x^{\tilde{u} - \epsilon-1} \bigl(\phi_1\bigl(\tilde{u} -\epsilon,
d(T_x)\bigr) - \phi_0\bigl(\tilde{u} - \epsilon,
d(T_x)\bigr)\bigr).
\end{eqnarray*}
Therefore, it suffices to show that the following converges to 0:
%
\begin{eqnarray}
\label{eqbigcor}
&&\Prob\Bigl( \sup_{u \leq\tilde{u} - \epsilon} x^{u-1}
\bigl(Z_1(uT_x) - \phi_1\bigl(u,
d(T_x)\bigr)\bigr)
\nonumber\\
&&\qquad{} + x^{\tilde{u} - \epsilon-1} \bigl(\phi_1\bigl(\tilde{u} -\epsilon ,
d(T_x)\bigr) - \phi_0\bigl(\tilde{u} - \epsilon,
d(T_x)\bigr)\bigr)
\\
&&\hspace*{18pt}\qquad{} + \sup_{u \leq\tilde{u} - \epsilon} x^{u-1}\bigl(\phi_0 \bigl( u,
d(T_x)\bigr) - Z_0^x(uT_x) \bigr)\nonumber
> 0 \Bigr).
\end{eqnarray}
To study this let us start by considering the first term in the sum above:
%
\begin{eqnarray}
\label{eqtemp1}
\qquad &&\sup_{u \leq\tilde{u} - \epsilon} \bigl|x^{u-1} \bigl(Z_1(uT_x)
- \phi_1\bigl(u, d(T_x)\bigr)\bigr)\bigr|
\nonumber
\\[-8pt]
\\[-8pt]
\nonumber
&&\qquad \leq \sup_{u \leq\tilde{u} - \epsilon} x^{u(1 + {\lambda_1}/{r}) -
\alpha} \sup_{u \leq\tilde{u} -
\epsilon} \bigl|
x^{\alpha- \lambda_1 u/r - 1} \bigl(Z_1(uT_x) - \phi_1
\bigl(u, d(T_x)\bigr)\bigr)\bigr|.
\nonumber
\end{eqnarray}
The second term in the product converges to zero in probability via
Theorem~\ref{thmapprox1}. The first term tends to zero by the
following argument:
\begin{eqnarray*}
&&\log \Bigl[\sup_{u \leq\tilde{u} - \epsilon} x^{u(1 + {\lambda
_1}/{r}) - \alpha} \Bigr] \\
&&\qquad= \log
\bigl[x^{(\tilde{u} - \epsilon)(1 +
{\lambda_1}/{r}) - \alpha} \bigr]
\\
&&\qquad = \log \biggl[ x^{-(1 + {\lambda_1}/{r}) \epsilon- \alpha} \exp \biggl[\tilde{u} \biggl(1 +
\frac{\lambda_1}{r}\biggr) \log x\biggr] \biggr]
\\
&&\qquad= \biggl(-\alpha- \biggl(1 + \frac{\lambda_1}{r}\biggr) \epsilon \biggr) \log{x}
+ \frac{\log x}{rT_x} \log \biggl[1 + x^\alpha \biggl(\frac{\lambda_1 +
r}{\mu}
\biggr) \biggr]
\\
&&\qquad\leq \biggl(-\alpha- \biggl(1 + \frac{\lambda_1}{r}\biggr) \epsilon \biggr) \log
{x} + \frac{\log x}{rT_x} \log \biggl[2x^\alpha \biggl(\frac{\lambda_1 +
r}{\mu}
\biggr) \biggr]
\\
&&\qquad= \alpha\frac{\log x}{rT_x} (\log x-r T_x ) - \biggl(
\frac
{\lambda_1}{r} + 1 \biggr)\epsilon\log x + \frac{\log x}{rT_x} \log \biggl[
\frac{2(\lambda_1 + 1)}{\mu} \biggr],
\end{eqnarray*}
where in the third equality we have utilized the fact that
\[
\tilde{u} \biggl(\frac{\lambda_1 + r}{r} \biggr) = \frac{1}{rT_x} \log \biggl[
\frac{ (\lambda_1 + r)x^\alpha}{\mu} + 1 \biggr]
\]
due to the definition of $\tilde{u}$. Observe that $ - (\frac
{\lambda_1}{r} + 1  )\epsilon\log x$
diverges to negative infinity, while the first and third terms approach
finite limits. This can be seen by observing that
\[
\frac{\log x}{rT_x} \rightarrow1
\]
in probability. Thus, we conclude that \eqref{eqtemp1} goes to zero in
probability.

The second term in \eqref{eqbigcor} is considered next. Via the
definition of $\phi_0(u,\break d(T_x))$ and the limit result on the
extinction time \eqref{eqexttimeconv}, we have that as $x\to\infty$,
\[
x^{\tilde{u} - \epsilon-1} \phi_0\bigl(\tilde{u} - \epsilon,
d(T_x)\bigr)\Rightarrow c^{-(\tilde u-\ve)}e^{-(\tilde u-\ve)\eta},
\]
where $\eta$ is a standard Gumbel random variable and $c$ is the
positive Yaglom constant. Importantly, this limit random variable is
positive with probability one.
The first term can be shown to approach zero by noting that
\[
x^{\tilde{u} -\epsilon-1} \phi_1\bigl(\tilde{u} - \epsilon,
d(T_x)\bigr) = x^{(\tilde{u} - \epsilon)( 1 + {\lambda_1}/{r} ) - \alpha} x^{\alpha- (\tilde{u} - \epsilon) \lambda_1 / r - 1}
\phi_1\bigl(\tilde {u}-\epsilon, d(T_x)\bigr),
\]
where the first term approaches zero, as argued previously since its
log approaches negative infinity, and the product of the remaining
terms approaches a constant times the exponential of a Gumbel, which is
again a result of \eqref{eqexttimeconv}.
The third term in \eqref{eqbigcor},
\[
\sup_{u \leq\tilde{u} - \epsilon} x^{u-1}\bigl(\phi_0 \bigl( u,
d(T_x)\bigr) - Z_0^x(uT_x)
\bigr),
\]
converges to zero in probability by Theorem~\ref{thmapprox1}. Therefore,
\[
\limsup_{x\to\infty}\Prob(\xi/T_x \leq\tilde{u} - \epsilon) \leq
\Prob \bigl( c^{-(\tilde u-\ve)}e^{-(\tilde u-\ve)\eta}\leq0\bigr) = 0.
\]

Next, we need to show that
$
\Prob(\xi/T_x \geq\tilde{u} + \epsilon)\rightarrow0.
$
We have by definition of $\xi$ that
\begin{eqnarray*}
&&\Prob(\xi/T_x > \tilde{u} + \epsilon) \leq\Prob\bigl(
Z_0\bigl((\tilde{u} + \epsilon)T_x \bigr) -
Z_1\bigl((\tilde{u} + \epsilon)T_x\bigr) > 0\bigr)
\\
& &\qquad= \Prob \bigl( x^{\alpha- (\tilde{u} + \epsilon)({\lambda_1}/{r}) -
1} \bigl(Z_0\bigl((\tilde{u} +
\epsilon)T_x \bigr) - Z_1\bigl((\tilde{u} + \epsilon
)T_x\bigr) \bigr) > 0\bigr)
\\
&&\qquad= \Prob\bigl( x^{\alpha- (\tilde{u} + \epsilon)({\lambda_1}/{r}) -
1} \bigl(Z_0\bigl((\tilde{u} +
\epsilon)T_x\bigr) - \phi_0\bigl(\tilde{u} + \epsilon,
d(T_x)\bigr)\bigr)
\\
&&\hspace*{10pt}\qquad\quad{} + x^{\alpha- (\tilde{u} + \epsilon)({\lambda_1}/{r}) -
1}\bigl(\phi_0\bigl(\tilde{u} + \epsilon,
d(T_x)\bigr) - \phi_1\bigl(\tilde{u} + \epsilon,
d(T_x)\bigr)\bigr)
\\
&&\hspace*{10pt}\qquad\quad{} + x^{\alpha- (\tilde{u} + \epsilon)({\lambda
_1}/{r}) - 1}\bigl(\phi_1\bigl(\tilde{u} + \epsilon,
d(T_x)\bigr) - Z_1\bigl((\tilde{u} +
\epsilon)T_x\bigr)\bigr) > 0 \bigr).
\end{eqnarray*}
It is easily shown that the right-hand side of the previous display
goes to 0 using analogous arguments from the analysis of \eqref{eqbigcor}.

\subsection{\texorpdfstring{Proof of Lemma \protect\ref{lemturnaroundtime1}}{Proof of Lemma 4}}
Here we establish\vspace*{1pt} Lemma~\ref{lemturnaroundtime1}, namely, that
$u^*(t)$ approximates $\tau/(\frac{1}{r}\log x+t)$.\vadjust{\goodbreak}

\begin{pf}
We will prove first that
%
\begin{equation}\label{eqp1}
\qquad \Prob \Bigl(\inf_{u\in[u^*(t)+\ve,\infty]} Z_1
\bigl(us_x(t)\bigr) < Z_0\bigl(u^*(t)s_x(t)
\bigr)+Z_1\bigl(u^*(t)s_x(t)\bigr) \Bigr) \rightarrow0.
\end{equation}
Consider the following decomposition of the event of interest,
\begin{eqnarray*}
&& \Prob \Bigl(\inf_{u\in[u^*(t)+\ve,\infty]} Z_1\bigl(us_x(t)
\bigr) < Z_0\bigl(u^*(t)s_x(t)\bigr)+Z_1
\bigl(u^*(t)s_x(t)\bigr) \Bigr)
\\
&&\qquad\leq \Prob \Bigl(\inf_{u\in[u^*(t)+\ve,\infty]}Z_1
\bigl(us_x(t)\bigr)<f\bigl(u^*(t)+\ve /2\bigr) \Bigr)
\\
&&\qquad\quad{} + \Prob \bigl(Z_0\bigl(u^*(t)s_x(t)
\bigr)+Z_1\bigl(u^*(t)s_x(t)\bigr)>f_{x,t}
\bigl(u^*(t)+\ve /2\bigr) \bigr).
\end{eqnarray*}
We can apply Markov's inequality to the last probability to see
\begin{eqnarray*}
&&\Prob \bigl(Z_0\bigl(u^*(t)s_x(t)\bigr)+Z_1
\bigl(u^*(t)s_x(t)\bigr)>f_{x,t}\bigl(u^*(t)+\ve /2\bigr)
\bigr)
\\
&&\qquad\leq\frac{\EE[Z_0(u^*(t)s_x(t))+Z_1(u^*(t)s_x(t))]}{f_{x,t}(u^*(t)+\ve
/2)}
\\
&&\qquad= \frac{f_{x,t}(u^*(t),t)}{f_{x,t}(u^*(t)+\ve/2)}
\\
&&\qquad= O\bigl(x^{-\lambda_1\ve/2r}\bigr),
\end{eqnarray*}
where the last equality follows from the ``steepness'' at the minimum
property~\eqref{eqsteepprop}.

Define the event
\[
A_{\ve}(x,t)=\inf\bigl\{Z_1\bigl(us_x(t)
\bigr)\dvtx u\in[u^*(t)+\ve,\infty)\bigr\}<f\bigl(u^*(t)+\ve/2\bigr).
\]
Then,
\begin{eqnarray}
\label{eqturnaroundineq1}
\Prob \bigl(A_{\ve}(x,t) \bigr)
&=&
\Prob \bigl(A_{\ve}(x,t), Z_1 \bigl(\bigl(u^*(t)+\ve
\bigr)s_x(t) \bigr)<f_{x,t}\bigl(u^*(t)+3\ve/4\bigr) \bigr)
\nonumber
\\
&&{} + \Prob \bigl(A_{\ve}(x,t), Z_1 \bigl(\bigl(u^*(t)+\ve
\bigr)s_x(t) \bigr)>f_{x,t}\bigl(u^*(t)+3\ve/4\bigr) \bigr)
\nonumber
\\[-8pt]
\\[-8pt]
\nonumber
&\leq& \Prob \bigl( Z_1 \bigl(\bigl(u^*(t)+\ve\bigr)s_x(t)
\bigr)<f_{x,t}\bigl(u^*(t)+3\ve /4\bigr) \bigr)
\\
&&{} + \Prob \bigl(A_{\ve}(x,t), Z_1 \bigl(
\bigl(u^*(t)+\ve\bigr)s_x(t) \bigr)>f_{x,t}\bigl(u^*(t)+3
\ve/4\bigr) \bigr).\nonumber
\end{eqnarray}
Using Chebyshev's inequality and the result in \eqref{eqsteepprop}
again, we see that
%
\begin{eqnarray}
\label{eqcheby1}
&&\Prob \bigl( Z_1 \bigl(\bigl(u^*(t)+\ve
\bigr)s_x(t) \bigr)<f_{x,t}\bigl(u^*(t)+3\ve /4\bigr) \bigr)
\nonumber
\\
&&\qquad= \Prob\bigl( Z_1 \bigl(\bigl(u^*(t)+\ve\bigr)s_x(t)
\bigr) - \EE\bigl[Z_1 \bigl(\bigl(u^*(t)+\ve \bigr)s_x(t)
\bigr)\bigr]
\nonumber
\\[-8pt]
\\[-8pt]
\nonumber
&&\hspace*{22pt}\qquad <f_{x,t}\bigl(u^*(t)+3\ve/4\bigr) - \EE\bigl[Z_1
\bigl(\bigl(u^*(t)+\ve\bigr)s_x(t) \bigr)\bigr] \bigr)
\\
&&\qquad\leq\frac{\EE[|Z_1 ((u^*(t)+\ve)s_x(t) ) - \EE Z_1
((u^*(t)+\ve)s_x(t) )|^2]}{|f_{x,t}(u^*(t) + {3 \epsilon
}/{4})-\EE Z_1 ((u^*(t)+\ve)s_x(t) )|^2}.\nonumber
\end{eqnarray}
Let us consider first the denominator in the above expression, and note
that since $u^*(t)$ minimizes $\EE Z_0(us_x(t)) + \EE Z_1(us_x(t))$, we
have that
\[
\frac{r}{\lambda_1} x e^{\lambda_0 u^*(t) (({1}/{r}) \log x + t)} \biggl( 1 - \frac{\mu}{x^\alpha(\lambda_1 + r)} \biggr) =
\frac
{x^{1-\alpha} \mu}{\lambda_1 + r}e^{\lambda_1 u^*(t)(({1}/{r}) \log x
+ t)}.
\]
Thus,
\begin{eqnarray*}
&&\EE\bigl[Z_1 \bigl(\bigl(u^*(t)+\ve\bigr)s_x(t)
\bigr)\bigr]
\\
&&\qquad = x e^{\lambda_0 u^*(t)(({1}/{r}) \log x + t)} \biggl[ -\frac{\mu
}{x^\alpha(\lambda_1 + r)}e^{\lambda_0 \ve(({1}/{r}) \log x + t)}
\\
&&\hspace*{106pt}\qquad{} + \frac{r}{\lambda_1} \biggl( 1 - \frac{\mu
}{x^\alpha(\lambda_1 + r)} \biggr)e^{\lambda_1 \ve(({1}/{r} )\log x
+ t)}
\biggr].
\end{eqnarray*}
Also,
\begin{eqnarray*}
&&f_{x,t}\bigl(u^*(t)+3 \ve/4\bigr)
\\
&&\qquad=xe^{\lambda_0u^*(t)(({1}/{r})\log x+t)} \biggl(1-\frac{\mu}{x^{\alpha
}(\lambda_1+r)} \biggr) \bigl(e^{\lambda_0 3 \ve/4(({1}/{r})\log
x+t)}\\
&&\hspace*{176pt}\qquad\quad{}+(r/
\lambda_1)e^{\lambda_1 3 \ve/4(({1}/{r})\log x +t)}\bigr),
\end{eqnarray*}
and therefore
\begin{eqnarray*}
&&\biggl|f_{x,t}\biggl(u^*(t) + \frac{3 \epsilon}{4}\biggr)-\EE
Z_1 \bigl(\bigl(u^*(t)+\ve \bigr)s_x(t)
\bigr)\biggr|^2
\\
&&\qquad= \biggl| x^{1-u^*(t)} e^{\lambda_0 u^*(t) t}\\
&&\qquad\quad{}\times \biggl[ \biggl(1-\frac{\mu
}{x^{\alpha}(\lambda_1+r)}
\biggr) \biggl( x^{-3\ve/4} e^{3\lambda_0 \ve
t/4} + \frac{r}{\lambda_1}
x^{3 \lambda_1 \ve/4r} e^{3\lambda_1 \ve
t/4} \biggr)
\\
&&\qquad\hspace*{25pt} { }+ \frac{\mu}{x^\alpha(\lambda_1 + r)}x^{-\ve} e^{-\ve t} -
\frac{r}{\lambda_1} \biggl( 1 - \frac{\mu}{x^\alpha(\lambda_1 + r)} \biggr)x^{\lambda_1 \ve/r}
e^{\lambda_1 \ve t} \biggr] \biggr|^2
\\
&&\qquad= \Omega\bigl(x^{2({\lambda_1 \ve}/{r} + 1 - u^*(t))}\bigr).
\end{eqnarray*}
Next we consider the variance term. For ease of notation define $\theta
(t) \equiv(u^*(t) + \ve)(\frac{1}{r} \log x + t)$, then from item
(iii) of Lemma~\ref{lemmamoments},
\begin{eqnarray*}
&&\operatorname{Var}\bigl[Z_1 \bigl(\bigl(u^*(t)+\ve
\bigr)s_x(t) \bigr)\bigr]
\\
&&\qquad= \biggl(\frac{\mu}{x^\alpha} \biggr)^2 \int_0^{\theta(t)}
\int_0^{\theta(t)}\operatorname{Cov}
\bigl(Z_0(s),Z_0(y)\bigr) e^{\lambda_1(2\theta(t)-s-y)} \,ds \,dy
\\
&&\qquad\quad{} + \biggl(\frac{\mu}{x^\alpha} \biggr) \int_0^{\theta(t)}
\EE Z_0(s) \bigl(2 r_1 e^{2 \lambda_1 (\theta(t) - s)} -
(r_1 + d_1) e^{\lambda_1 (\theta(t) - s)} \bigr) \,ds.
\end{eqnarray*}
We now establish that $\operatorname{Cov}(Z_0(s),Z_0(y))$ is an $O(x)$ quantity. Since
the population $Z_0$ starts with $x$ independent cells, we can write
the covariance as
\begin{eqnarray*}
\operatorname{Cov}\bigl(Z_0(s),Z_0(y)\bigr)&=&\EE \Biggl[
\sum_{j=1}^xZ_0^{(j)}(y)
\sum_{i=1}^xZ_0^{(i)}(s)
\Biggr]-\EE\bigl[Z_0(s)\bigr]\EE\bigl[Z_0(y)\bigr]
\\
&=& x\EE_1\bigl[Z_0(s)Z_0(y)\bigr]+x(x-1)
\EE_1\bigl[Z_0(s)\bigr]\EE_1
\bigl[Z_0(y)\bigr]\\
&&{}-\EE\bigl[Z_0(s)\bigr]\EE
\bigl[Z_0(y)\bigr]
\\
&=& x\operatorname{Cov}_1\bigl(Z_0(s),Z_0(y)
\bigr)=O(x).
\end{eqnarray*}
Based on this we know that the second term in the definition of
$\operatorname{Var}Z_1$ is the dominant term and therefore $\operatorname{Var}[Z_1
((u^*(t)+\ve)s_x(t) )]=O (x^{1-\alpha+2 \lambda_1
{(u^*(t)+\ve)}/{r}} )$.

Then, in order to establish that \eqref{eqcheby1} goes to zero it
suffices to show that
\[
\biggl(1-\alpha+2 \lambda_1\frac{u^*(t)+\ve}{r}-2\biggl(
\frac{\lambda_1 \ve
}{r} + 1 - u^*(t)\biggr) \biggr)\log x \rightarrow-\infty.
\]
The result in the previous display follows from the definition of
$u^*(t)$ in \eqref{eqturnaround}, and therefore we can conclude that
\[
\lim_{x\to\infty}\Prob \bigl( Z_1 \bigl(\bigl(u^*(t)+\ve
\bigr)s_x(t) \bigr)<f_{x,t}\bigl(u^*(t)+3\ve/4\bigr)
\bigr)=0.
\]

It remains to show that the final probability in display \eqref
{eqturnaroundineq1} is negligible for large $x$.
Observe that if we start out with a collection of $n$ independent
cells, each following branching processes with net-growth rate $\lambda_1>0$,
then by the law of large numbers the fraction, $f_n$, of those
cells whose lineage eventually dies out satisfies the following limit:
$\lim_{n\to\infty}f_n = p_E(\lambda_1)<1$, where $p_E(\lambda_1)$ is
the probability of a single cell's descendants going extinct and is
strictly less than 1 because $\lambda_1>0$.
Therefore, define $\rho_x(u,t)$ to be the fraction of type-1 cells
present at time $u(\frac{1}{r}\log x+t)$ that eventually die out.
Notice then that in order for the event described in the last line of
display \eqref{eqturnaroundineq1} to occur it is necessary that
\[
\rho_x\bigl(u^*(t)+\ve,t\bigr) \geq1-\frac{f_{x,t}(u^*(t)+\ve
/2)}{f_{x,t}(u^*(t)+3\ve/4)}.
\]
Then from the ``steepness'' property we have that, for $x$ sufficiently
large, $\rho_x(u^*(t)+\ve,t)>p_E(\lambda_1)+\eta$, for some $\eta>0$.
Of course, from the law of large numbers we have that
\begin{eqnarray*}
\Prob \bigl( \rho_x\bigl(u^*(t)+\ve,t\bigr)>p_E(
\lambda_1)+\eta, Z_1 \bigl(\bigl(u^*(t)+\ve
\bigr)s_x(t) \bigr)>f_{x,t}\bigl(u^*(t)+3\ve/4\bigr) \bigr)
\end{eqnarray*}
converges to $0$ as $x\to\infty$,
thus establishing \eqref{eqp1}.

Moving on we next establish that
%
\begin{equation}
\label{eqp2} \qquad\Prob \Bigl(\inf_{u\in[0,u^*(t)-\ve]} Z_0
\bigl(us_x(t)\bigr) < Z_0\bigl(u^*(t)s_x(t)
\bigr)+Z_1\bigl(u^*(t)s_x(t)\bigr) \Bigr) \rightarrow0
\end{equation}
as $x \rightarrow\infty$. Note that based on arguments from the above
case it suffices to establish that the following probability converges
to 0 as $x\to\infty$:
\begin{eqnarray*}
&&\Prob \Bigl(\inf_{u\in[0,u^*(t)-\ve]}Z_0\bigl(us_x(t)
\bigr)<f_{x,t}\bigl(u^*(t)-\ve /2\bigr) \Bigr)
\\
&&\qquad\leq
\frac{1}{(c_0-1)f_{x,t}(u^*(t)-\ve/2)} \\
&&\qquad\quad{}\times\EE\bigl[ \bigl(c_0f_{x,t}\bigl(u^*(t)-
\ve/2\bigr)-Z_0 \bigl(\bigl(u^*(t)-\ve\bigr)s_x(t) \bigr)
\bigr)^+ \bigr]
\\
&&\qquad\leq \frac{c_0f_{x,t}(u^*(t)-\ve/2)}{(c_0-1)f_{x,t}(u^*(t)-\ve/2)}\\
&&\qquad\quad{}\times \Prob \bigl(Z_0 \bigl(\bigl(u^*(t)-\ve
\bigr)s_x(t) \bigr)<c_0f_{x,t}\bigl(u^*(t)-\ve/2
\bigr) \bigr),
\end{eqnarray*}
where we chose $c_0>1$, and the penultimate inequality follows from
Doob's inequality and that $-Z_0(\cdot s_x(t))$ is a submartingale. The
final probability can be rewritten as
\begin{eqnarray*}
&&\Prob\bigl(Z_0 \bigl(\bigl(u^*(t)-\ve\bigr)s_x(t) \bigr)-\EE
Z_0 \bigl(\bigl(u^*(t)-\ve \bigr)s_x(t) \bigr)
\\
&&\qquad < c_0f_{x,t}\bigl(u^*(t)-\ve/2\bigr)-\EE Z_0
\bigl(\bigl(u^*(t)-\ve\bigr)s_x(t) \bigr)\bigr). 
\end{eqnarray*}
Note since $\EE Z_0 ((u^*(t)-\ve)s_x(t) ) = x^{1+\ve
-u^*(t)}e^{\lambda_0 t(u^*(t)-\ve)}$ and $f_{x,t}(u^*(t)-\ve/2) =
x^{1+\ve/2-u^*(t)}(1+o(1))$, we know that there exists a positive
constant $C_0$ such that for $x$ sufficiently large
\[
\EE Z_0 \bigl(\bigl(u^*(t)-\ve\bigr)s_x(t)
\bigr)-cf_{x,t}\bigl(u^*(t)-\ve/2\bigr)\geq C_0x^{1+\ve-u^*(t)},
\]
and therefore for $x$ sufficiently large
\begin{eqnarray*}
&&\Prob\bigl(Z_0 \bigl(\bigl(u^*(t)-\ve\bigr)s_x(t) \bigr)-\EE
Z_0 \bigl(\bigl(u^*(t)-\ve \bigr)s_x(t) \bigr)
\\
&&\quad < cf_{x,t}\bigl(u^*(t)-\ve/2\bigr)-\EE Z_0 \bigl(
\bigl(u^*(t)-\ve\bigr)s_x(t) \bigr)\bigr)
\\
&&\qquad\leq \Prob \bigl(\bigl\llvert Z_0 \bigl(\bigl(u^*(t)-\ve
\bigr)s_x(t) \bigr)-\EE Z_0 \bigl(\bigl(u^*(t)-\ve
\bigr)s_x(t) \bigr)\bigr\rrvert >C_0x^{1+\ve-u^*(t)}
\bigr).
\end{eqnarray*}
Thus, via Chebyshev's inequality we have
\begin{eqnarray*}
&&\Prob \bigl(\bigl\llvert Z_0 \bigl(\bigl(u^*(t)-\ve
\bigr)s_x(t) \bigr)-\EE Z_0 \bigl(\bigl(u^*(t)-\ve
\bigr)s_x(t) \bigr)\bigr\rrvert >C_0x^{1+\ve-u^*(t)}
\bigr)
\\
&&\qquad\leq \frac{x^{1/2} (\operatorname{Var}_1Z_0 ((u^*(t)-\ve)s_x(t) )
)^{1/2}}{C_0x^{1+\ve-u^*(t)}}
\\
&&\qquad= O \bigl(x^{-(1+\ve-u*)/2} \bigr),
\end{eqnarray*}
where the final equality follows by evaluating the variance of
$Z_0(u^*(t)s_x(t))$.~%
\end{pf}

\subsection{\texorpdfstring{Proof of Lemma \protect\ref{lemturnaroundunif}}{Proof of Lemma 5}}
As in the proof of the previous lemma, we consider the deviations to
the left and right of $u^*(t)$ separately.\vadjust{\goodbreak} First, note that if $t_0 >
t_1$, then
\begin{eqnarray*}
&&\inf_{u \in[0, u^*(t_0) - \ve)} Z_0\biggl(u\biggl(\frac{1}{r}\log
x+t_0\biggr)\biggr) + Z_1\biggl(u\biggl(
\frac{1}{r}\log x+ t_0\biggr)\biggr)
\\
&&\qquad= \inf\biggl\{Z_0(s) + Z_1(s) \dvtx s \leq
\bigl(u^*(t_0) - \ve\bigr) \biggl(\frac{1}{r} \log x +
t_0 \biggr) \biggr\}
\\
&&\qquad \geq \inf\biggl\{Z_0(s) + Z_1(s) \dvtx s \leq
\bigl(u^*(t_1) - \ve\bigr) \biggl(\frac{1}{r} \log x +
t_1 \biggr) \biggr\}
\\
&&\qquad= \inf_{u \in[0, u^*(t_1) - \ve)}Z_0\biggl(u\biggl(\frac{1}{r}\log x+
t_1\biggr)\biggr) + Z_1\biggl(u\biggl(
\frac{1}{r}\log x+ t_1\biggr)\biggr).
\end{eqnarray*}
Furthermore, it follows from the definition of $u^*(t)$ that
$Z_i(u^*(t)(\frac{1}{r}\log x+t))=Z_i(u^*(s)(\frac{1}{r}\log x+s))$ for
all $s,t$. In particular, $Z_i (u^*(s)(\frac{1}{r}\log x+s)
)=Z_i(u^*T_x)$.
Then via Lemma \eqref{lemturnaroundtime1},
\[
\Prob \Bigl(\sup_{t \in[-M,M]} \inf_{u \in[0, u^*(t) - \ve
)}Z_0
\bigl(us_x(t)\bigr) + Z_1\bigl(us_x(t)\bigr)
< Z_0\bigl(u^*T_x\bigr) + Z_1
\bigl(u^*T_x\bigr) \Bigr)
\]
converges to zero as $x \rightarrow\infty$. A similar argument can be
used for deviations to the right of $u^*(t)$ to show that
\[
\Prob \Bigl(\sup_{t \in[-M,M]} \inf_{u \in[u^*(t) + \ve, \infty
)}Z_0
\bigl(us_x(t)\bigr) + Z_1\bigl(us_x(t)\bigr)
<Z_0\bigl(u^*T_x\bigr) + Z_1
\bigl(u^*T_x\bigr) \Bigr)
\]
also converges to zero as $x \rightarrow\infty$, establishing the lemma.


%


\printaddresses

\end{document}